\documentclass[a4paper,10pt]{article}
 \usepackage{benstyle}
 
\title{Stable reconstructions in Hilbert spaces and the resolution of the Gibbs phenomenon}
\author{Ben Adcock \\ Department of Mathematics \\ Simon Fraser University \\ Burnaby, BC V5A 1S6 \\ Canada \\ \hspace{20pc} \and Anders C. Hansen \\ DAMTP, Centre for Mathematical Sciences \\ University of Cambridge \\ Wilberforce Rd, Cambridge CB3 0WA \\ United Kingdom}
\begin{document}
\maketitle

\begin{abstract}
We introduce a method to reconstruct an element of a Hilbert space in terms of an arbitrary finite collection of linearly independent reconstruction vectors, given a finite number of its samples with respect to any Riesz basis.  As we establish, provided the dimension of the reconstruction space is chosen suitably in relation to the number of samples, this procedure can be numerically implemented in a stable manner.  Moreover, the accuracy of the resulting approximation is completely determined by the choice of reconstruction basis, meaning that the reconstruction vectors can be tailored to the particular problem at hand.

An important example of this approach is the accurate recovery of a piecewise analytic function from its first few Fourier coefficients.  Whilst the standard Fourier projection suffers from the Gibbs phenomenon, by reconstructing in a piecewise polynomial basis, we obtain an approximation with root exponential accuracy in terms of the number of Fourier samples and exponential accuracy in terms of the degree of the reconstruction function.  Numerical examples illustrate the advantage of this approach over other existing methods.
\end{abstract}

\section{Introduction}
Suppose that $\rH$ is a separable Hilbert space with inner product $\left < \cdot , \cdot \right >$ and corresponding norm $\nm{\cdot}$.  In this paper, we consider following problem: given the first $m$ samples $\{\left < f , \psi_j \right >\}_{j=1}^m$ of an element $f \in \rH$ with respect to some Riesz basis $\{\psi_j\}^{\infty}_{j=1}$ of $\rH$ (the \textit{sampling} basis), reconstruct $f$ to high accuracy.  Not only does such a problem lie at the heart of modern sampling theory \cite{eldar2005general,unser2000sampling}, it also occurs in a myriad of applications, including image processing (in particular, Magnetic Resonance Imaging), and the numerical solution of hyperbolic partial differential equations (PDEs).

In practice, straightforward reconstruction of $f$ may be achieved via orthogonal projection with respect to the sampling basis.  Indeed, for an arbitrary $f \in \rH$, this is the best possible strategy.  However, in many important circumstances, this approximation converges only slowly in $m$, when measured in the norm on $\rH$, or not at all, if a stronger norm (for example, the uniform norm) is considered.

A prominent instance of this problem is the recovery of a function $f:\bbT \rightarrow \bbR$ from its first $m$ Fourier coefficients (here $\bbT= [-1,1)$ is the unit torus).  In this instance, $\rH = \rL^{2}(\bbT)$ is the space of all $2$-periodic, square-integrable functions of one variable.  Provided $f$ is analytic, it is well-known that its Fourier series (the orthogonal projection with respect to the Fourier basis) converges exponentially fast.  However, whenever $f$ has a jump discontinuity,  its Fourier series  suffers from the well-known Gibbs phenomenon \cite{korner}.  Whilst convergence occurs in the $\rL^{2}$ norm, uniform convergence is lacking, and the approximation is polluted by characteristic $\ord{1}$ oscillations near the discontinuity.  Moreover, the rate of convergence is also slow: only $\ordu{m^{-\frac{1}{2}}}$ when measured in the $\rL^2$ norm, and $\ord{m^{-1}}$ pointwise away from the discontinuity.  Needless to say, the Gibbs phenomenon is a significant blight of many practical applications of Fourier series \cite{jerri1998gibbs}.  It is a testament to its importance that the design of effective techniques for its removal remains an active area of inquiry \cite{GottGibbsRev,Tadmor1}.

Returning to the general form of the problem, let us now suppose that some additional information is known about the function $f$.  For example, $f$ may be sparse in a particular basis (e.g. a polynomial of low degree) or may possess certain regularity.  In particular, in the Fourier instance, we may know that $f$ is piecewise analytic with jump discontinuities at known locations in $\bbT$.  In this circumstance, it seems plausible that a better approximation to $f$ can be obtained by expanding in a different basis (e.g. a piecewise polynomial basis).  To this end, let us introduce the so-called \textit{reconstruction} space (of dimension $n$) and seek to approximate $f$ by an element $f_{n,m}$ consisting of $n$ linearly independent elements of this space.  

As we will show in due course, provided reconstruction is carried out in a certain manner, a suitable approximation $f_{n,m}$ can always be found.  Essential to this approach is that $m$ (the number of samples) is chosen sufficiently large in comparison to $n$ (equivalently, $n$ is chosen sufficiently small in comparison to $m$).  However, provided this is the case, the approximation $f_{n,m}$ inherits the principal features of the reconstruction space.  In particular, $f_{n,m}$ is quasi-optimal (or, under certain conditions, asymptotically optimal), in sense that the error $\| f - f_{n,m} \|$ can be bounded by a constant multiple of $\| f - \cQ_{n} f \|$, where $\cQ_{n} f$ is the orthogonal projection onto the reconstruction space (in other words, the best approximation to $f$ from this space).  Moreover, from a practical standpoint, this method can be implemented by solving a linear least squares problem.  Whenever the reconstruction vectors are suitably chosen (e.g. if they form a Riesz basis), the corresponding linear system is well-conditioned and the least squares problem can be solved in $\ord{m n}$ operations by standard iterative techniques.

Consider once more the example of Fourier series.  Let $f:\bbT \rightarrow \bbR $ be an analytic function  with jump discontinuity at $x=-1$ (equivalently, an analytic, nonperiodic function).  As mentioned, the Fourier series of $f$ lacks uniform convergence.  However, since $f$ is analytic, it makes sense to seek to reconstruct $f$ in a system of polynomials.  It is well-known that the best $n^{\rth}$ degree polynomial approximation of an analytic function converges exponentially fast in $n$ \cite{boyd}.  As we shall prove, with $n = \ord{\sqrt{m}}$, we obtain a quasi-optimal polynomial approximation $f_{n,m}$ to $f$ from only its first $m$ Fourier coefficients, regardless of the particular family of polynomials used.  This results in exponential convergence of $f_{n,m}$ to $f$ in $n$ (the polynomial degree), or root exponential convergence in $m$ (the number of Fourier samples).    Moreover, whenever Legendre polynomials are used, the approximation is computable in a stable manner in only $\ord{m n}$ operations.  The use of, for example, Chebyshev polynomials results in a method with $\ord{n}$ condition number, that can be implemented in $\cO(m n^{\frac{3}{2}})$ operations.  Furthermore, whilst retaining the aforementioned features, this procedure can be easily generalised to recovery of a piecewise analytic function of one variable (using piecewise polynomial bases), and to the case of multivariate functions defined in tensor-product regions.

There are a number of existing algorithms for the removal of the Gibbs phenomenon from Fourier series.  One of the most well-known, which also provides a polynomial approximant, is the Gegenbauer reconstruction technique \cite{GottGibbs3,GottGibbsRev,GottGibbs1}.  As we discuss further in Section \ref{Fourrecovsec}, the method developed in this paper has a number of key advantages over this device.  Numerical results also indicate its superior performance for a number of test problems.

The method developed in this paper was previously introduced by the authors in \cite{BAACHShannon} within the context of abstract sampling theory.  Whilst this problem, in this abstract form, has been extensively studied in the last couple of decades (in particular, by Eldar et al \cite{eldar2003sampling,eldar2005general}, see also \cite{unser2000sampling}), to the best of our knowledge this particular method does not appear in any existing literature.  For a more detailed discussion of the relation of this approach to existing schemes we refer the reader to \cite{BAACHShannon}.  Conversely, in this paper, after presenting the general version of the method in abstract terms, we will focus primarily on its application to the Fourier coefficient reconstruction problem.  On this topic, a similar approach, but only dealing with reconstructions in Legendre polynomials  from Fourier samples of analytic functions, was discussed in \cite{hrycakIPRM}.  This can be viewed as a special case of our general framework.  Furthermore, by examining this example as part of this framework, we are able to extend and improve the work of \cite{hrycakIPRM} in the following ways: (i) we derive a procedure allowing for reconstructions in any polynomial basis, not just Legendre polynomials, (ii) we extend this approach to reconstructions of piecewise smooth functions using (arbitrary) piecewise polynomial bases, (iii) we generalise this work to smooth functions of arbitrary numbers of variables and (iv) we obtain improved estimates for both the error and the necessary scaling $n = \ord{\sqrt{m}}$ required for implementation.

Aside from yielding these improvements, a great benefit of the general framework presented in this paper is that it is immediately applicable to a whole host of other reconstruction problems.  To illustrate this generality, in the final part of this paper we consider its application to the accurate reconstruction of a piecewise analytic function from its orthogonal polynomial expansion coefficients.  Such a problem is typical of that occurring in the application of polynomial spectral methods to hyperbolic PDEs \cite{gottlieb2001spectral}, where the shock formation inhibits fast convergence of the polynomial approximation.  As we highlight, this issue can be overcome in a completely stable fashion by reconstructing in a piecewise polynomial basis.

The outline of the remainder of this paper is as follows.  In Section \ref{genrecsect} we introduce the reconstruction procedure and establish both stability and error estimates.  Section \ref{Fourrecovsec} is devoted to (piecewise) polynomial reconstructions from Fourier samples.  In Section \ref{higherdimsect} we consider reconstructions from tensor-product spaces, and in Section \ref{othersampsect} we discuss other reconstruction problems.  Finally, in Section \ref{conclusions} we present open problems and challenges.

\section{General theory of reconstruction}
\label{genrecsect}
In this section, we describe the reconstruction procedure in its full generality.  To this end, suppose that $\{ \psi_j \}^{\infty}_{j=1}$ is a Riesz basis (the sampling basis) for a separable Hilbert space $\rH$ over the field $\bbC$.  Let $\left < \cdot , \cdot \right >$ be the inner product on $\rH$, with associated norm $\nm{\cdot}$.  Recall that, by definition, $\spn \{ \psi_1 , \psi_{2} , \ldots \}$ is dense in $\rH$ and 
\be{
\label{Rieszprop}
c_1 \sum^{\infty}_{j=1} | \alpha_{j} |^2 \leq \left \| \sum^{\infty}_{j=1} \alpha_{j} \psi_{j} \right \|^2 \leq c_2 \sum^{\infty}_{j=1} | \alpha_{j} |^2,\quad \forall \alpha = \{\alpha_1,\alpha_{2},\ldots \} \in l^2(\bbN),
}  
for positive constants $c_1,c_2$.  Equivalently, $\psi_{j} = \cU (\Psi_{j})$, where $\{ \Psi_{j} \}^{\infty}_{j=1}$ is an orthonormal basis for $\rH$ and $\cU : \rH \rightarrow \rH$ is a bounded, bijective operator.  Using this definition, it is easy to deduce that $\{ \psi_{j} \}^{\infty}_{j=1}$ also satisfies the frame property
\be{
\label{frameprop}
d_1 \| f \|^2 \leq \sum^{\infty}_{j=1} | \left < f , \psi_{j} \right > |^2 \leq d_2 \| f \|^2,\quad \forall f \in \rH,
} 
for $d_1,d_2>0$, where the largest possible value for $d_2$ is $\| \cU \|^2_{\rH \rightarrow \rH}$ and the smallest possible value for $d_1$ is $\| \cU^{-1} \|^{-2}_{\rH \rightarrow \rH}$ \cite{christensen2003introduction}.

Suppose now that the first $m$ coefficients of an element $f \in \rH$ with respect to the sampling basis are given:
\be{
\label{fsamples}
\hat{f}_{j} = \left < f , \psi_{j} \right >,\quad j=1,\ldots,m.
}
Set $\rS_{m} = \spn \{ \psi_1,\ldots,\psi_m \}$ and let $\cP_{m} : \rH \rightarrow \rS_{m}$ be the mapping
\be{
\label{Pmdef}
f \mapsto \cP_{m} f = \sum^{m}_{j=1} \left < f , \psi_{j} \right > \psi_{j}.
}  
We now seek to reconstruct $f$ in a different basis.  To this end, suppose that $\{ \phi_{1},\ldots,\phi_{n} \}$ are linearly independent reconstruction vectors and define $\rT_{n} = \spn \{ \phi_1,\ldots,\phi_n \}$.  Let $\cQ_{n} : \rH \rightarrow \rT_{n}$ be the orthogonal projection onto $\rT_{n}$.  Direct computation of $\cQ_{n} f$, the best approximation to $f$ from $\rT_{n}$, is not possible, since the coefficients $\left < f , \phi_j \right >$ are unknown.  Instead, we seek to use the values \R{fsamples} to compute an approximation $f_{n,m} \in \rT_{n}$ that is \textit{quasi-optimal}, i.e. $\| f - \cQ_{n} f \| \leq \| f - f_{n,m} \| \leq C \| f - \cQ_{n} f \|$ for some constant $C>0$ independent of $f$, $n$ and $m$.  
To do this, we introduce the sesquilinear form $a_m : \rH \times \rH \rightarrow \bbC$, given by
\be{
\label{sesq}
a_{m}(g,h) = \left < \cP_{m} g ,  h \right >,\quad \forall g,h \in \rH.
}
Note that, since
\bes{
\left < \cP_{m} g ,  h \right > = \sum^{m}_{j=1} \left < f , \psi_j \right > \overline{\left < g , \psi_{j} \right >} = \overline{a_m(g,f)},\quad \forall f,g \in \rH,
}
$a_{m}$ is a Hermitian form on $\rH \times \rH$ (here $\overline{z}$ is the complex conjugate of $z \in \bbC$).  With this to hand, we now define $f_{n,m}$ by the following condition
\be{
\label{recoveqn1}
a_{m}(f_{n,m},\phi) = a_m(f,\phi),\quad \forall \phi \in \rT_{n}.
}
Upon setting $\phi = \phi_j$, $j=1,\ldots,n$, this becomes an $n \times n$ linear system of equations for the coefficients $\alpha_1,\ldots,\alpha_n$ of $f_{n,m} = \sum^{n}_{j=1} \alpha_j \phi_j$.  We shall defer a discussion of the computation of this approximation to Section \ref{comp_f}: in the remainder of this section we consider the analysis of $f_{n,m}$.

Before proving the main theorem regarding \R{recoveqn1}, let us first give an explanation as to why this approach works.  As mentioned, key to this technique is that the parameter $m$ is sufficiently large in comparison to $n$.  To this end, let $n$ be fixed and suppose that $m \rightarrow \infty$.  Due to \R{Rieszprop}, the mappings $\cP_{m}$ converge strongly to a bounded, linear operator $\cP$, given by
\be{
\label{Pdef}
\cP f = \sum^{\infty}_{j=1} \left < f , \psi_{j} \right > \psi_j,\quad \forall f \in \rH.
}
Hence, for large $m$, the equations \R{recoveqn1} defining $f_{n,m}$ resemble the equations 
\be{
\label{tilde}
a(\tilde f_{n} , \phi) = a (f , \phi), \quad \forall \phi \in \rT_n,
}
where $a : \rH \times \rH \rightarrow \bbC$ is the Hermitian form $a(f,g) = \langle \cP f , g \rangle$.  Thus, it is reasonable to expect that $f_{n,m} \rightarrow \tilde f_n$ as $m \rightarrow \infty$, 
provided such a function $\tilde f_n$ exists.  However,
\thm{
For all $n \in \bbN$, the function $\tilde{f}_{n}$ exists and is unique.  Moreover,
\be{
\label{fterror}
\| f - \tilde{f}_{n} \| \leq \frac{d_2}{d_1} \| f - \cQ_{n} f \|,
}
where $d_1$ and $d_2$ arise from (\ref{frameprop}).
}
This theorem can be established with a straightforward application of the Lax--Milgram theorem and its counterpart, C\'ea's lemma \cite{gilbarg2001elliptic}.  Indeed, due to \R{frameprop} and \R{Pdef},
\be{
\label{IPequiv}
d_{1} \| g\|^2 \leq a(g,g) = \sum^{\infty}_{j=1} | \left < g , \psi_{j} \right > |^2 \leq d_2 \| g \|^2,\quad \forall g \in \rH.
}
Hence the form $a(\cdot,\cdot)$ defines an equivalent inner product on $\rH$.  Nonetheless, we shall present a self-contained proof, since similar techniques will be used subsequently.
\prf{
Let $\cU : \rT_{n} \rightarrow \bbC^{n}$ be the linear mapping $g \mapsto  \{ \langle \cP g , \phi_j \rangle \}^{n}_{j=1}$.  To prove existence and uniqueness of $\tilde{f}_{n}$ it suffices to show that $\cU$ is invertible, upon which it follows that $\tilde f_{n} = \cU^{-1} \{ \langle \cP f , \phi_j \rangle \}^{n}_{j=1} $.  Suppose that $\cU g = 0$.  Then, by definition, $\langle \cP g , \phi_j \rangle = 0$ for $j=1,\ldots,n$.  Using linearity, we deduce that $\langle \cP g , g \rangle = 0$.  Now, it follows from \R{IPequiv} that $0 = \langle  \cP g , g \rangle  \geq d_1 \| g \|^2$, giving $g=0$.  Hence, $\cU$ is invertible and $\tilde f_n$ exists and is unique.

Now consider the error estimate \R{fterror}.  Using \R{frameprop} once more, we obtain
\bes{
\| f - \tilde{f}_{n} \|^2 \leq \frac{1}{d_1} \left < \cP(f-\tilde{f}_{n}) , f  - \tilde{f}_{n} \right >.
}
By definition of $\tilde{f}_{n}$, $\langle \cP(f-\tilde{f}_{n}) , \phi \rangle = 0$, $\forall \phi \in \rT_{n}$.  In particular, setting $\phi = \tilde{f}_{n} - \cQ_{n} f$, yields
\bes{
\| f - \tilde{f}_{n} \|^2 \leq \frac{1}{d_1} \left < \cP(f-\tilde{f}_{n}) , f - \cQ_{n} f \right >.
}
Since $a(\cdot , \cdot)$ forms an equivalent inner product on $\rH$, an application of the Cauchy--Schwarz inequality gives
\bes{
\| f - \tilde{f}_{n} \|^2 \leq \frac{1}{d_1}\left [ \langle  \cP (f-\tilde{f}_{n}) , f - \tilde{f}_{n} \rangle 
\langle \cP(f-\cQ_{n} f ) , f-\cQ_{n} f \rangle\right ]^{\frac{1}{2}} \leq \frac{d_{2}}{d_1} \| f - \tilde{f}_{n} \| \| f - \cQ_{n} f \|,
}
as required.
}
This theorem establishes existence and quasi-optimality of $\tilde{f}_{n} \approx f_{n,m}$, thereby giving an intuitive argument for the success of this method.  We now wish to fully confirm this observation.  To this end, let
\be{
\label{Cdef}
C_{n,m} = \inf_{\substack{\phi \in \rT_{n} \\ \| \phi \| =1}} \sum^{m}_{j=1} | \left < \phi , \psi_{j} \right > |^2 =  \inf_{\substack{\phi \in \rT_{n} \\ \| \phi \| =1}}  \left < \cP_{m} \phi , \phi \right > = \inf_{\substack{\phi \in \rT_{n} \\ \| \phi \| =1}} a_{m}(\phi,\phi).
} 
The quantity $C_{n,m}$ plays a fundamental role in this paper.  Note that
\lem{
\label{Cproplem}
For all $n,m \in \bbN$, $0 \leq C_{n,m} \leq d_2$.  Moreover, for each $n$, $C_{n,m} \rightarrow C^{*}_{n} \geq d_1$ as $m \rightarrow \infty$, where
\bes{
C^{*}_{n} = \inf_{\substack{\phi \in \rT_{n} \\ \| \phi \| =1}} \sum^{\infty}_{j=1} | \left < \phi , \psi_{j} \right > |^2 =  \inf_{\substack{\phi \in \rT_{n} \\ \| \phi \| =1}}  \left < \cP \phi , \phi \right > = \inf_{\substack{\phi \in \rT_{n} \\ \| \phi \| =1}} a(\phi,\phi),
}
and $d_1$ is defined in (\ref{frameprop}).
}
\prf{
Consider first the quantity
\be{
\label{epdef}
\epsilon_{n,m} =  \sup_{\substack{\phi \in \rT_{n} \\ \| \phi \| =1}} \sum_{j > m} | \left < \phi , \psi_{j} \right >|^2 = \sup_{\substack{\phi \in \rT_{n} \\ \| \phi \| =1}}  \left < \cP \phi - \cP_{m} \phi , \phi \right >.
}
Due to \R{frameprop}, the infinite sum is finite for any fixed $\phi$, and tends to zero as $m \rightarrow \infty$.  Now let $\{ \Phi_j \}^{n}_{j=1}$ be an orthonormal basis for $\rT_{n}$ and set $\phi = \sum^{n}_{j=1} \alpha_{j} \Phi_{j}$.  Two applications of the Cauchy--Schwarz inequality gives
\bes{
\sum_{j > m} | \left < \phi , \psi_{j} \right >|^2 \leq \| \phi \|^2 \sum^{n}_{k=1} \sum_{j>m} | \left < \Phi_{k} , \psi_{j} \right > |^2 .
}
Hence $\epsilon_{n,m} \leq \sum^{n}_{k=1} \sum_{j>m} | \left < \Phi_{k} , \psi_{j} \right > |^2 $, and we deduce that $\epsilon_{n,m}$ is both finite and $\epsilon_{n,m} \rightarrow 0 $ as $m \rightarrow \infty$.  Noticing that $| C_{n,m} - C^{*}_{n} | \leq \epsilon_{n,m}$, $\forall n,m \in \bbN$, gives the first part of the proof.  For the second, we merely use \R{frameprop}.
}
Aside from $C_{n,m}$, we also define the quantity $D_{n,m}$ by
\be{
\label{Ddef}
D_{n,m} =\sup_{\substack{f \in \rT^{\perp}_{n} \\ \| f \| =1 }} \sup_{\substack{g \in \rT_{n} \\ \| g \|=1 }}  \left | \left < \cP_{m} f , g \right > \right |.
}
For this, we have the following lemma:
\lem{
\label{dlem}
For all $m,n \in \bbN$, $0\leq D_{n,m} \leq d_{2}$.  Moreover, suppose that $\cP$ is such that $\cP(\rT_{n}) \subseteq \rT_{n}$ (for example, when $\cP = \cI$ is the identity), then $D^2_{n,m} \leq c_{2} \epsilon_{n,m}$, where $\epsilon_{n,m}$ is as in \R{epdef}.  In particular, for fixed $n$, $D_{n,m} \rightarrow 0$ as $m \rightarrow \infty$.
}
\prf{
Let $f,g \in \rH$.  By definition,
\be{
\label{contstep}
\left < \cP_{m} f , g \right > = \sum^{m}_{j=1} \left < f , \psi_{j} \right > \overline{\left < g , \psi_{j} \right >} \leq \left [ \sum^{m}_{j=1} | \left < f , \psi_{j} \right > |^2 \right ]^{\frac{1}{2}} \left [ \sum^{m}_{j=1} | \left < g, \psi_{j} \right > |^2 \right ]^{\frac{1}{2}}.
}
Hence \R{frameprop} now gives the first result.  Now suppose that $f \in \rT^{\perp}_{n}$ and $\cP(\rT_{n}) \subseteq \rT_{n}$.  Since $\cP_{m}$ is self-adjoint, we have $\left < \cP_{m} f , g \right > = \left < f , \cP_{m} g \right >  = \left < f , \cP_{m} g - \cP g \right >$.  Here the second equality is due to the fact that $f \perp \rT_{n}$ and $\cP g \in \rT_{n}$ for $g \in \rT_{n}$.  By the Cauchy--Schwarz inequality, we obtain
\bes{
D_{n,m} \leq \sup_{\substack{g \in \rT_{n} \\ \| g \| =1}}\| \cP g - \cP_{m} g \|.
}
For $g \in \rT_{n}$, we note from \R{Rieszprop} that
\bes{
\| \cP g - \cP_{m} g \|^2 \leq c_{2} \sum_{j > m} | \left < g , \psi_{j} \right > |^2  = c_{2} \left < \cP g - \cP_{m} g , g \right > \leq c_{2} \epsilon_{n,m} \| g \|^2,
}
where the final equality follows from the definition of $\epsilon_{n,m}$.
}
At this moment, we mention the following important fact.  The constants $C_{n,m}$, $D_{n,m}$ and $C^{*}_{n}$, as well as the approximation $f_{n,m}$, are determined only by the space $\rT_{n}$, not by the choice of reconstructions vectors $\phi_{1},\ldots,\phi_{n}$ themselves.  As we shall discuss later, the choice of reconstruction vectors only affects the stability of the scheme.

This observation aside, we are now able to state the main theorem of this section:
\thm{
\label{reconthm} 
For every $n \in \bbN$ there exists an $m_0$ such that the approximation $f_{n,m}$, defined by \R{recoveqn1}, exists and is unique for all $m \geq m_0$, and satisfies the stability estimate $\| f_{n,m} \| \leq d_2 C^{-1}_{n,m} \| f\|$.  Furthermore,
\be{
\label{reconerrest}
\| f - \cQ_{n} f \| \leq \| f - f_{n,m} \| \leq K_{n,m} \| f - \cQ_n f \|,\quad K_{n,m} = \sqrt{1+D^{2}_{n,m} C^{-2}_{n,m}}.
}
Specifically, the parameter $m_0$ is the least value of $m$ such that $C_{n,m}>0$.
}
To prove this theorem, we first recall that a Hermitian form $a : \rH \times \rH \rightarrow \bbR$  is said to be \textit{continuous} if, for some constant $\gamma>0$, $|a(f,g)| \leq \gamma \| f \| \| g \|$ for all $f,g \in \rH$.  Moreover, $a$ is \textit{coercive}, provided $a(f,f) \geq \omega \| f \|^2$, $\forall f \in \rH$, for $\omega>0$ constant \cite{gilbarg2001elliptic}.  We now require the following lemma:
\lem{
\label{amcoerclem} 
Suppose that $a_{m} : \rH \times \rH \rightarrow \bbR$ is the sesquilinear form $a_{m}(f,g) = \left < \cP_m f , g \right >$.  Then, $a_{m}$ is continuous with constant $\gamma \leq d_{2}$.  Moreover, for every $n \in \bbN$ there exists an $m_0$ such that the restriction of $a_{m}$ to $\rT_{n} \times \rT_{n}$ is coercive for all $m \geq m_0$.  Specifically, if $C_{n,m}$ is given by \R{Cdef}, then $m_0$ is the least value of $m$ such that $C_{n,m} >0$, and, for all $m \geq m_0$, $a_{m}(f,f) \geq C_{n,m} \| f \|^2$, $\forall f \in \rT_{n}$.  Finally, for all $f \in \rH$ and $g \in \rT_{n}$, we have $a_{m}(f-\cQ_{n} f , g )  \leq D_{n,m} \| f - \cQ_{n} f \|  \| g \|$.
}
\prf{
Continuity follows immediately from \R{contstep}.  For the second and final results, we merely use the definitions \R{Cdef} and \R{Ddef} of $C_{n,m}$ and $D_{n,m}$ respectively.
}
\prf{[Proof of Theorem \ref{reconthm}]
To establish existence and uniqueness, it suffices to prove that the linear operator $\cU : \rT_{n} \rightarrow \bbC^{n}$, $g \mapsto \{ \left < \cP_{m} g, \phi_{j} \right > \}^{n}_{j=1}$ is invertible.  Suppose that $g \in \rT_{n}$ with $\cU g = 0$.  By definition, we have $\left < \cP_{m} g , \phi_{j} \right > = 0$ for $j=1,\ldots,n$.  Using linearity, it follows that $\left < \cP_{m} g , g \right > = 0$.  Lemma \ref{amcoerclem} now gives $0 \leq C_{n,m} \| g \|^2 \leq 0$.  Hence $g = 0$, and therefore $\cU$ has trivial kernel.

Stability of $f_{n,m}$ is easily established from the continuity and coercivity conditions.  Setting $\phi = f_{n,m}$ in \R{recoveqn1} gives
\bes{
C_{n,m} \| f_{n,m} \|^2 \leq a_{m}(f_{n,m},f_{n,m})  = a_{m}(f,f_{n,m}) \leq d_2 \| f \| \| f_{n,m} \|,
}
as required.  Now consider the error estimate \R{reconerrest}.  Suppose that we define $e_{n,m} = f_{n,m} - \cQ_n f \in \rT_{n}$.  Then, by definition of $f_{n,m}$, we have $a_m(e_{n,m},\phi) = a_{m}(f-\cQ_n f , \phi)$, $\forall \phi \in \rT_{n}$.  In particular, setting $\phi = e_{n,m}$, we obtain
\be{
\label{enbound}
\| e_{n,m} \| \leq C^{-1}_{n,m} D_{n,m} \| f - \cQ_n f \|.
}
Since $\cQ_{n} f$ is the orthogonal projection onto $\rT_{n}$, we have $\| f - f_{n,m} \|^2 = \| e_{n,m} \|^2 + \| f - \cQ_{n} f \|^2$, which gives the full result.
}
When $a_m$ is shown to be continuous and coercive, it may be tempting to seek to apply the Lax--Milgram theorem and C\'ea's lemma  to obtain Theorem \ref{reconthm}.  However, the Hermitian form $a_m$, when considered as a mapping $\rH \times \rH \rightarrow \bbC$, will not, in general, be coercive.  This is readily seen from the definition of $C_{n,m}$.  The finite-dimensional operator $\cP_{m}|_{\rT_n}$ converges uniformly to $\cP |_{\rT_{n}}$, whereas its infinite-dimensional counterpart $\cP_{m}$ typically does not (for example, when $\cP_{m}$ is the Fourier projection operator and $\rH = \rL^{2}(\bbT)$).  Hence, $a_m$ only becomes coercive when restricted to $\rT_{n} \times \rT_{n}$, and these standard results do not automatically apply.

Though Theorem \ref{reconthm} establishes an estimate for the error $f-f_{n,m}$ measured in the natural norm on $\rH$, it is also useful to derive a result valid for any other norm defined on a suitable subspace of $\rH$ (for example, this may be the uniform norm on $\bbT$ in the case of Fourier series).  To this end, let $\tnm{\cdot}$ be such a norm and define $\rG = \{ g \in \rH : \tnm{g} < \infty \}$.  We have
\cor{
\label{unifcor}
Suppose that $f \in \rG$,  $\rT_{n} \subseteq \rG$ and that $f_{n,m}$ is defined by \R{recoveqn1}.  Then, for all $m \geq m_{0}$,
\be{
\label{uniferr}
\tnm{f-f_{n,m}} \leq \tnm{f-\cQ_{n} f } +\frac{k_n D_{n,m}}{C_{n,m}} \| f - \cQ_n f \|,
}
where $k_n = \sup_{\substack{\phi \in \rT_{n} \\ \| \phi \| =1 }} \tnm{\phi}$ and $C_{n,m}$, $D_{n,m}$ are given by \R{Cdef} and \R{Ddef} respectively.
} 
\prf{
Let $e_{n,m} = f_{n,m} - \cQ_n f$ once more.  Since $e_{n,m} \in \rT_{n}$, it follows from the definition of $k_n$ and the inequality \R{enbound} that
\bes{
\tnm{e_{n,m}} \leq k_n \| e_{n,m} \| \leq \frac{k_n D_{n,m}}{C_{n,m}} \| f - \cQ_n f \|.
}
The full result is obtained  from the triangle inequality $\tnm{f-f_{n,m}} \leq \tnm{e_{n,m}} + \tnm{ f - \cQ_n f }$.
}
\begin{remark}
In practice, it is useful to have an upper bound for the constant $k_{n}$.  A simple exercise gives $k_{n} \leq \sum^{n}_{j=1} \tnm{\Phi_{j}}$, where $\{ \Phi_j \}^{n}_{j=1}$ is any orthonormal basis for $\rT_{n}$.
\end{remark}

Theorem \ref{reconthm} confirms that the approximation $f_{n,m}$ is quasi-optimal whenever $m$ is sufficiently large in comparison to $n$.  In other words, the error $\| f - f_{n,m} \|$ can be bounded by a constant multiple (independent of $n$) of $\| f - \cQ_{n} f \|$.  Naturally, whenever $\{ \phi_{1},\ldots , \phi_{n} \}$ are the first $n$ vectors in an infinite sequence $\phi_{1},\phi_{2},\ldots$ that form a basis for $\rH$ (a natural assumption to make from a practical standpoint), we find that $f_{n,m}$ converges to $f$ at the same rate as $\cQ_{n}f$, provided $m$ scales appropriately with $n$.  Under the same assumptions, Corollary \ref{unifcor} also verifies convergence of $f_{n,m}$ to $f$ in $\tnm{\cdot}$, whenever $\cQ_{n} f \rightarrow f$ in this norm and $k_{n} \| f - \cQ_{n} f \| \rightarrow 0$ as $n \rightarrow \infty$.

Note on other feature of this framework.  Provided $m$ is such that $C_{n,m}>0$, any function $f \in \rT_{n}$ will be recovered exactly by $f_{n,m}$.  In the language of sampling theory, this property is commonly referred to as \textit{perfect} reconstruction.

\begin{remark}
Whilst this framework relies on the fact that $m$ can range independently of $n$, a natural question to ask is what happens if $m$ is set equal to $n$ (note that, in this case, one recovers the well-known \textit{consistent sampling} framework of Eldar et al \cite{eldar2005general,unser2000sampling}).  This question was discussed in detail in \cite{BAACHShannon}, where it demonstrated that such an approach often leads to severe ill-conditioning as $n=m \rightarrow \infty$.  Additionally, stringent restrictions are placed on the types of vectors $f$ that can be reconstructed (see also Section \ref{othermethsec}).  Conversely, by allowing $m$ to vary independently of $n$, we obtain a reconstruction $f_{n,m}$ that is guaranteed to converge for any vector $f \in \rH$.  Moreover, as we discuss in Section \ref{comp_f}, provided the reconstruction vectors are suitably chosen, the computation of $f_{n,m}$ is completely stable.
\end{remark}

\subsection{The case of an orthogonal sampling basis}
\label{orthsampsect}
In the previous setup, the sampling basis was assumed to be a Riesz basis.  The particular case of $\{ \psi_{j} \}^{\infty}_{j=1}$ being an orthonormal basis warrants further study, since this situation often arises in applications (e.g. Fourier sampling).  

In this setting, both \R{Rieszprop} and \R{frameprop} (which are now identical) hold with equality and with all constants being precisely one.  In other words, Parseval's identity
\bes{
\| f \|^2 = \sum^{\infty}_{j=1} | \left < f , \psi_{j} \right > |^2,\quad \forall f \in \rH,
}
holds.  Our first result provides a geometric interpretation for the constant $C_{n,m}$ in terms of subspace angles.  Recall that if $\rU$ and $\rV$ are subspaces of $\rH$, then the angle $\theta_{\rU \rV}$ is given by
\be{
\label{subang}
\cos (\theta_{\rU \rV}) = \inf_{\substack{u \in \rU \\ \| u \| = 1}} \| \cQ_{\rV} u \|, 
}
where $\cQ_{\rV} : \rH \rightarrow \rV$ is the orthogonal projection.  We have

\lem{
\label{subanglem}
Suppose that $\{ \psi_{j} \}^{\infty}_{j=1}$ is an orthonormal basis of $\rH$.  Then $C_{n,m} = \cos^2 (\theta)$, where $\theta = \theta_{\rT_{n} \rS_{m}}$ is the angle between the subspaces $\rT_{n}$ and $\rS_{m}$.
}
\prf{
Whenever $\{ \psi_{j} \}^{\infty}_{j=1}$ is orthonormal, the operator $\cP_{m} $ is the orthogonal projection onto $\rS_{m}$.  Hence the result follows immediately from \R{Cdef} and \R{subang}.
}
Whilst this property is interesting, the following result has much greater significance:
\thm{
\label{reconorththm}
Suppose that $\{ \psi_{j} \}^{\infty}_{j=1}$ is an orthonormal basis of $\rH$.  Then, for all $m \geq m_{0}$, the approximation $f_{n,m}$ satisfies
\be{
\label{reconerrestorth}
\| f - \cQ_{n} f \| \leq \| f - f_{n,m} \| \leq K_{n,m} \| f - \cQ_n f \|
}
where
\bes{
K_{n,m} = \sqrt{1+(1-C_{n,m})C^{-2}_{n,m}} = \sqrt{1+\tan^2 \theta \sec^2 \theta},
}
and $\theta$ is as in Lemma \ref{subanglem}.  In particular, for fixed $n$,  $K_{n,m} \rightarrow 1$ as $m \rightarrow \infty$, and hence $f_{n,m} \rightarrow \cQ_{n} f$ as $m \rightarrow \infty$.
}
\prf{
Note that $\cP = \cI$ for an orthonormal sampling basis.  In particular, $\epsilon_{n,m} = C^{*}_{n} - C_{n,m}$ and $C^{*}_{n} = 1$.  Hence, using Lemma \ref{dlem} we find that $D^2_{n,m} \leq 1 - C_{n,m}$.  The result now follows from Theorem \ref{reconthm}.
}

From this we conclude the following: not only is $f_{n,m}$ quasi-optimal, it is also \textit{asymptotically optimal} in the sense that $f_{n,m} \rightarrow \cQ_{n} f$, the best approximation to $f$ from $\rT_{n}$, as $m \rightarrow \infty$.  Hence, using this approach, we can recover an approximation to $f$ that is arbitrarily close to the error minimising approximant (which, as mentioned, cannot be computed directly from the given samples).  Moreover, the rate of convergence of $f_{n,m}$ to $\cQ_{n}f$ is completely independent of the particular vector $f$, and relies only on rate of decay of the parameter $1-C_{n,m}$.

Note that asymptotic optimality also occurs for general Riesz bases whenever $\cP(\rT_{n}) \subseteq \rT_{n}$, in which case $D_{n,m} \rightarrow 0$ as $m \rightarrow \infty$ (see Lemma \ref{dlem}).  The case of orthonormal sampling vectors presents the most obvious example of a basis satisfying this condition.

\begin{remark}
Whenever the vectors $\{ \psi_{j} \}^{\infty}_{j=1}$ are not orthonormal, a natural question to ask is whether we can modify the approach for computing $f_{n,m}$ to recover asymptotic optimality.  This can be easily done, at least in theory, by replacing the operator $\cP_{m}$, given by \R{Pmdef}, with the orthogonal projection $\rH \rightarrow \rT_{n}$.  In this case, both Lemma \ref{subanglem} and Theorem \ref{reconorththm} will hold for nonorthogonal sampling bases.  The downside of the approach is that it requires additional computational cost to compute $f_{n,m}$, as we explain at the end of the next section.

Another potential means to recover asymptotic optimality is to define $\cP_{m} g = \sum^{m}_{j=1} \left < g, \psi_{j} \right > \psi^{*}_{j}$, where $\{ \psi^{*}_{j} \}$ is the set of dual vectors to the sampling vectors $\{ \psi_{j} \}$.  In this case, $\cP_{m} \rightarrow \cI$ strongly, and asymptotic optimality follows.  In practice, however, one may not have access to the dual vectors, thus this approach cannot necessarily be readily implemented.
\end{remark}

\subsection{Oblique asymptotic optimality}
Whenever the sampling basis is orthonormal, $f_{n,m}$ converges to the best approximation $\cQ_{n} f$ as $m \rightarrow \infty$.  An obvious question to ask is what can be said in the general case, where the vectors $\{ \psi_{j} \} $ only form a Riesz basis?  The intuitive explanation given previously indicated that $f_{n,m} \approx \tilde{f}_{n}$, where $\tilde{f}_{n}$ is defined by \R{tilde}.  We now wish to confirm this observation.  In fact, as in the orthonormal case, we may demonstrate a stronger result: namely, for fixed $n \in \bbN$, $f_{n,m} \rightarrow \tilde{f}_{n}$ as $m\rightarrow \infty$ at a rate that is independent of the particular vector $f \in \rH$.  

Recall that the form $a(\cdot,\cdot)$ yields an equivalent inner product on $\rH$.  Since $\tilde{f}_{n}$ is defined by the equations $a(\tilde{f}_n , \phi ) = a(f,\phi)$, $\forall \phi \in \rT_{n}$, the mapping $f \mapsto \tilde{f}_{n}$ is the orthogonal projection onto $\rT_{n}$ with respect to this inner product.  Letting $\| g \|_{a} = \sqrt{a(g,g)}$ be the corresponding norm on $\rH$, we now define the constants
\be{
\label{Ctdef}
\tilde C_{n,m} = \inf_{\substack{\phi \in \rT_{n} \\ \| \phi \|_a=1}} \left < \cP_{m} \phi , \phi \right >,\qquad  \tilde D_{n,m} = \sup_{\substack{f \in \rT^{\perp}_n \\ \| f \|_a=1}} \sup_{\substack{g \in \rT_{n} \\ \| g \|_a=1}} | \left < \cP_{m}f , g \right > |.
}
In this instance, $\rT^{\perp}_{n}$ is defined with respect to the $a$-inner product, i.e. $\rT^{\perp}_{n} = \{ f \in \rH : a(f,\phi) = 0,\  \forall \phi \in \rT_{n} \}$.  Conversely, when considered with respect to the canonical inner product, this subspace is precisely $\cP(\rT_{n})^{\perp} = \{ f \in \cH : \langle f , \phi \rangle =0,\  \forall \phi \in \cP(\rT_{n}) \}$.

Note the similarity between $\tilde C_{n,m}$ and $\tilde{D}_{n,m}$ and the quantities $C_{n,m}$ and $D_{n,m}$ defined in \R{Cdef} and \R{Ddef} respectively.  Roughly speaking, the former measure the deviation of $f_{n,m}$ from $\cQ_{n} f$, whereas, as we will subsequently show, the latter determine the deviation of $f_{n,m}$ from $\tilde{f}_{n}$.  

With these definitions to hand, identical arguments to those given in the proofs of Lemmas \ref{Cproplem} and \ref{dlem} now yield:

\lem{
For all $m,n \in \bbN$, $\tilde{C}_{n,m} \geq \frac{1}{d_2} C_{n,m}$, where $C_{n,m}$ is as in \R{Cdef}.  Moreover, for fixed $n$, $\tilde{C}_{n,m} \rightarrow 1$ as $m \rightarrow \infty$.
}
\lem{
For all $m,n \in \bbN$, $\tilde D_{n,m} \leq d_2$ and $\tilde{D}^2_{n,m} \leq c_{2}(1-\tilde{C}_{n,m})$.  In particular, for fixed $n$, $\tilde D_{n,m} \rightarrow 0$ as $m \rightarrow \infty$.
}
Using these lemmas, we deduce
\cor{
If $f_{n,m}$ and $\tilde f_{n}$ are given by \R{recoveqn1} and \R{tilde} respectively, then 
\bes{
\| f_{n,m} - \tilde{f}_{n} \|_{a} \leq \frac{\tilde D_{n,m}}{\tilde C_{n,m}} \| f - \tilde{f}_{n} \|_{a}, 
}
and we have the error estimate
\bes{
\| f - \tilde f_{n} \|_{a} \leq \| f - f_{n,m} \|_{a} \leq \tilde{K}_{n,m} \| f - \tilde f_{n} \|_{a},\quad \tilde K_{n,m} = \sqrt{1+\tilde D^2_{n,m}\tilde C^{-2}_{n,m} }.
}
In particular, for any $f \in \rH$, $f_{n,m} \rightarrow \tilde{f}_{n}$ as $m \rightarrow \infty$.
}
\prf{
Since $(f_{n,m} - \tilde f_{n}) \in \rT_{n}$, we have
\bes{
\tilde C_{n,m} \| f_{n,m} - \tilde f_{n} \|^2_{a} \leq   \left < \cP_{m} (f_{n,m} - \tilde{f}_{n}) , f_{n,m}-\tilde{f}_{n} \right >.
}
Moreover, because $\langle \cP_{m} f_{n,m} , \phi \rangle = \langle \cP_{m} f , \phi \rangle $, we deduce that
\bes{
\tilde C_{n,m} \| f_{n,m} - \tilde f_{n} \|^2_{a} \leq  \left < \cP_{m} ( f - \tilde{f}_{n} ) , f_{n,m} - \tilde{f}_{n} \right > \leq \tilde D_{n,m} \| f - \tilde{f}_{n} \|_a \| f_{n,m} - \tilde{f}_{n} \|_a,
}
where the second inequality follows from the definition \R{Ctdef} of $\tilde D_{n,m}$ and the fact that $(f-\tilde{f}_{n}) \in \rT^{\perp}_{n}$, the orthogonal complement of $\rT_{n}$ with respect to the $a$-inner product.
}
Note that the mapping $\cW_{n}: f \mapsto \tilde{f}_{n}$ is an \textit{oblique} projection with respect to the inner product $\langle \cdot , \cdot \rangle$ on $\rH$.  In particular, $\cW_{n}$ has range $\rT_{n}$ and kernel $\cP(\rT_n)^{\perp}$, and we have the decomposition $\rH = \rT_{n} \oplus \cP(\rT_{n})^{\perp}$.  For this reason, we say that $f_{n,m}$ possesses \textit{oblique asymptotic} optimality.

\subsection{Computation of $f_{n,m}$}
\label{comp_f}
Recall that the computation of the approximation $f_{n,m}$ involves solving the system of equations \R{recoveqn1}.  This can be interpreted as the normal equations of a least squares problem.  Suppose that $f_{n,m} = \sum^{n}_{j=1} \alpha_j \phi_j$, $\alpha = (\alpha_1,\ldots,\alpha_n) \in \bbC^n$ and $\hat{f} = (\hat{f}_{1},\ldots,\hat{f}_{m})$.  If $U$ is the $m \times n$ matrix with $(j,k)^{\rth}$ entry $\left < \phi_k , \psi_{j} \right >$, then \R{recoveqn1} is given exactly by $A \alpha = U^{\dag} \hat{f}$, where $A = U^{\dag} U $ and $U^{\dag}$ is the adjoint of $U$.  In other words, the vector $\alpha$ is the least squares solution of the problem $U \alpha \approx \hat{f}$.

This system can be solved iteratively by applying conjugate gradient iterations to the normal equations, for example.  The number of required iterations is dependent on the condition number $\kappa(A)$ of the matrix $A$.  Specifically, the number of iterations required to obtain numerical convergence (i.e. to within a prescribed tolerance) is proportional to $\sqrt{\kappa(A)}$ \cite{golub}.  In particular, if $\kappa(A)$ is $\ord{1}$ for all $n$ and $m \geq m_0$, then the number of iterations is also $\ord{1}$ for all $n$.  Hence, the cost of computing $f_{n,m}$ is determined solely by the number of operations required to perform matrix-vector multiplications involving $U$.  In other words, only $\ord{m n}$ operations.  

Naturally, aside from this consideration, the condition number of $A$ is also important since it determines susceptibility of the numerical computation to both round-off error and noise.  Specifically, an error of magnitude $\epsilon$ in the inputs (i.e. the samples $\hat{f}_{j}$, $j=1,\ldots,m$) will yield an error of magnitude roughly $\kappa(A) \epsilon$ in the output $f_{n,m}$.

For these reasons it is of utmost importance to study the condition number of $A$.  For this, we first introduce the Hermitian matrix $\tilde{A} \in \bbC^{n \times n}$ with $(j,k)^{\rth}$ entry $\left < \phi_j , \phi_k \right >$.  Note that $\tilde{A}$ is the Gram matrix of the vectors $\{ \phi_1 , \ldots, \phi_n \}$.  In particular, $\kappa(\tilde{A})$ is a measure of the suitability of the particular vectors in which to compute $\cQ_n f$.  With $\tilde A$ to hand, we also introduce the related matrix $\tilde{A}_a \in \bbC^{n \times n}$ with $(j,k)^{\rth}$ entry $a(\phi_j,\phi_k) = \left < \cP \phi_j , \phi_k \right >$, i.e. the Gram matrix with respect to the inner product $a(\cdot,\cdot)$.

The following lemma comes as no surprise:

\lem{
The matrices $\tilde{A}$ and $\tilde{A}_a$ are spectrally equivalent.  In particular, for all $n \in \bbN$,
\bes{
\frac{d_1}{d_2} \kappa(\tilde{A}) \leq \kappa (\tilde A_a) \leq \frac{d_2}{d_1} \kappa (\tilde{A}).
}
}
\prf{
For any Hermitian matrix $B$, the condition number is the ratio of the largest and smallest eigenvalues (in absolute value).  Moreover, if $B$ is positive definite, then 
\be{
\label{Atevals}
\inf_{\substack{\alpha \in \bbC^n \\ \alpha \neq 0}} \left \{\frac{\alpha^{\dag} B \alpha}{\alpha^\dag \alpha} \right \} = \lambda_{\min}(B),\quad \sup_{\substack{\alpha \in \bbC^n \\ \alpha \neq 0}} \left \{\frac{\alpha^{\dag} B \alpha}{\alpha^\dag \alpha} \right \} = \lambda_{\max}(B).
}
If $\phi =  \sum^{n}_{j=1} \alpha_j \phi_j$, then $\alpha^{\dag} \tilde A \alpha = \| \phi \|^2$ and $\alpha^{\dag} \tilde{A}_a \alpha = a(\phi,\phi)$.  Hence, spectral equivalence now follows immediately from \R{IPequiv}.
}
Concerning the condition number of the matrix $A$, we now have the following:
\lem{
\label{condnumlem}
Suppose that $m \geq m_0$, where $m_0$ is as in Theorem \ref{reconthm}, and $\tilde{C}_{n,m}$ and $C_{n,m}$ are given by \R{Cdef} and \R{Ctdef} respectively.  Then
\bes{
\tilde{C}_{n,m} \kappa(\tilde{A}_a) \leq \kappa(A) \leq \frac{1}{\tilde C_{n,m}} \kappa(\tilde{A}_a),\qquad \frac{C_{n,m}}{d_2} \kappa(\tilde{A}) \leq \kappa(A) \leq \frac{d_2}{ C_{n,m}} \kappa(\tilde A).
}
Moreover, for fixed $n$, $A \rightarrow \tilde{A}_a$ as $m \rightarrow \infty$, and, if $\cP = \cI$, $A \rightarrow \tilde{A} = \tilde{A}_{a}$.
}
\prf{
The matrix $A$ is Hermitian and, provided $m \geq m_{0}$, positive definite.  Hence, its eigenvalues are given by \R{Atevals}.  For $\phi = \sum^{n}_{j=1} \alpha_{j} \phi_{j}$, we have $\alpha^{\dag} A \alpha = \langle \cP_{m} \phi , \phi \rangle$.  By definition of $C_{n,m}$ and $\tilde{C}_{n,m}$, we find that $\lambda_{\min}(A) \geq C_{n,m} \lambda_{\min}(\tilde{A})$ and $ \lambda_{\min}(A) \geq \tilde C_{n,m} \lambda_{\min}(\tilde{A}_a)$.  Moreover, by \R{Rieszprop} we have $\lambda_{\max}(A) \leq d_2 \lambda_{\max}(\tilde{A})$ and $\lambda_{\max}(A) \leq \lambda_{\max}(\tilde{A}_a)$.  The first result now follows immediately from \R{Atevals}.  For the second, we merely note that each entry of $A$ converges to the corresponding entry of $\tilde A_{a}$ as $m \rightarrow \infty$.
}
Note the important conclusion of this lemma: computing $f_{n,m}$ from \R{recoveqn1} is no more ill-conditioned than the computation of the orthogonal projection $\cQ_{n}f$ or the oblique projection $\cW_{n} f$ in terms of the vectors $\{ \phi_1 , \ldots,\phi_{n} \}$.  In practice, it is often true that these vectors correspond to the first $n$ vectors in a basis $\{ \phi_j \}^{\infty}_{j=1}$ of $\rH$ with additional structure.  Whenever this is the case, as the following trivial corollary indicates, we can expect good conditioning:

\cor{
\label{reconreiszstab}
Suppose that $\{ \phi_j \}^{\infty}_{j=1}$ is a Riesz basis for $\rH$ (with respect to $\langle \cdot , \cdot \rangle$) with constants $c'_1$ and $c'_{2}$.  Then 
\bes{
\kappa(A) \leq \frac{c'_2 d_2}{c'_1C_{n,m}}.
}
}
\prf{
This follows immediately follows from \R{Rieszprop} and Lemma \ref{condnumlem}.
}
Put together, a fundamental conclusion of Theorem \ref{reconthm}, Lemma \ref{condnumlem} and Corollary \ref{reconreiszstab} is the following: for a given reconstruction space $\rT_{n}$, the individual vectors $\phi_{1},\ldots,\phi_{n}$ can be chosen arbitrarily, without altering the analysis of the approximation $f_{n,m}$ (which itself does not depend on the individual vectors used to represent it).  The choice of vectors only becomes important when considering the condition number of linear system to solve.  Moreover, the quality of a system of vectors for the reconstruction problem is completely intrinsic, in that it is determined only by the corresponding Gram matrix (in particular, it is independent of the sampling space).

Corollary \ref{reconreiszstab} confirms that the approximation $f_{n,m}$ can be readily computed in a stable manner for many choices of reconstruction basis.  However, to fully implement this method, as we discuss further in the next section, it is useful to have numerical way of computing $C_{n,m}$.  The following lemma provides such a means:
\lem{
\label{Cformlem}
The quantity $C_{n,m}$ is given by $C_{n,m} =\lambda_{\min}(\tilde{A}^{-1} A)$.  Moreover, if $\tilde{A}$ and $A$ commute, then $C_{n,m} =  1 - \| I - \tilde{A}^{-1} A \|$.  In particular, if $\{ \phi_j \}^{n}_{j=1}$ is an orthonormal basis, then $C_{n,m} = \lambda_{\min}(A) = 1 - \| I - A \|$.
}
\prf{
By definition $C_{n,m} = \inf_{\substack{\phi \in \rT_{n} \\ \| \phi \| =1}} \left < \cP_{m} \phi , \phi \right > $.  Letting $\phi = \sum^{n}_{j=1} \alpha_{j} \phi_{j}$, we find that
\bes{
C_{n,m} = \inf_{\substack{\alpha \in \bbC^n \\ \alpha \neq 0}} \frac{\sum^{n}_{j,k=1} \alpha_{j} \overline{\alpha}_{k}  \left < \cP_{m} \phi_{j} , \phi_{k} \right > }{\sum^{n}_{j,k=1} \alpha_{j} \overline{\alpha}_{k}  \left < \phi_{j} , \phi_{k} \right > } =  \inf_{\substack{\alpha \in \bbC^n \\ \alpha \neq 0}} \frac{\alpha^{\dag}  A \alpha}{\alpha^{\dag}  \tilde{A} \alpha}.
}
We now claim that, for arbitrary Hermitian positive definite matrices $B$ and $C$ with $B$ nonsingular, the following holds:
\bes{
\inf_{\substack{\alpha \in \bbC^n \\ \alpha \neq 0}} \frac{\alpha^{\dag} C \alpha}{\alpha^{\dag} B \alpha} = \lambda_{\min}(B^{-1} C), \quad \sup_{\substack{\alpha \in \bbC^n \\ \alpha \neq 0}} \frac{\alpha^{\dag} C \alpha}{\alpha^{\dag} B \alpha} = \lambda_{\max}(B^{-1} C).
}
To do so, write $B = D^{\dag} D$, with $D$ nonsingular.  Then, after rearranging, we obtain
\bes{
\inf_{\substack{\alpha \in \bbC^n \\ \alpha \neq 0}} \frac{\alpha^{\dag} C \alpha}{\alpha^{\dag} B \alpha} = \inf_{\substack{\beta \in \bbC^n \\ \beta \neq 0}} \frac{\beta^{\dag} D^{-\dag} C D^{-1} \beta}{\beta^\dag \beta} = \lambda_{\min}( D^{-\dag} C D^{-1} ),
}
for example.  However, a trivial calculation confirms that the eigenvalues of $D^{-\dag} C D^{-1}$ are identical to those of $B^{-1} C$, thus establishing the claim.  Since $\tilde{A}$ is nonsingular, this confirms that $C_{n,m} = \lambda_{\min}(\tilde{A}^{-1} A)$.  For the second result, we merely notice that $\lambda_{\min}(B) = 1 - \lambda_{\max}(I-B) = 1 - \| I - B \|$, whenever $B$ is Hermitian.
}

In Section \ref{orthsampsect}, we briefly discussed a modified approach where the operator $\cP_{m}$, usually given by \R{Pmdef}, was replaced by the orthogonal projection operator.  The advantage of this approach is that it guarantees asymptotic optimality.  However, the downside is additional computational expense.  Indeed, the corresponding matrix is of the form $A = U^{\dag} V^{-1} U$, where $V \in \bbC^{m \times m}$ has $(j,k)^{\rth}$ entry $\left < \psi_{j} , \psi_{k} \right >$.  Hence, if conjugate gradients iterations are used, at each stage we are required to compute matrix-vector products involving the $m \times m$ matrix $V^{-1}$ (assuming that $V^{-1}$ had been precomputed).  In general, this requires $\ord{m^2}$ operations.  Thus, we incur a cost of $\ord{m^2}$, as opposed to $\ord{m n}$ for the original algorithm.  Hence, in practice it may be better settle for only quasi- and oblique asymptotic optimality, whilst retaining a lower computational cost.

\subsection{Conditions for guaranteed, quasi-optimal recovery}
Let us return to the standard form of the method once more.  To implement this method, it is necessary to have conditions that guarantee nonsingularity, stability and quasi-optimal recovery.  In other words, for given sampling and reconstruction bases, we wish to study the quantity
\be{
\label{Theta}
\Theta(n; \theta) = \min \left \{ m \in \bbN : C_{n,m} \geq \theta \right \},\quad \theta \in (0,d_2),
}
where $C_{n,m}$ is given by \R{Cdef} and $d_2$ stems from \R{frameprop} .  Note that $C_{n,m} \leq d_2$ by \R{frameprop}, thereby explaining the stated range of $\theta$.  Also, by Lemma 
\ref{Cproplem}, we have that $\lim_{m\rightarrow \infty} C_{n,m} \geq d_1>0$, thus $\Theta$ is well-defined.  

By definition, $\Theta(n;\theta)$ is the least $m$ such that $\| f- f_{n,m} \| \leq c(\theta) \| f - \cQ_{n} f \|$, 
where 
\be{
\label{smallc} 
c(\theta) = \sqrt{1+d^2_2 \theta^{-2}} \quad \text{or} \quad \sqrt{1+(1-\theta)\theta^{-2}},
} 
whenever the sampling basis is orthonormal.  In other words, the least $m$ required for quasi-optimal recovery with constant $c(\theta)$.  Thus, provided $m \geq \Theta(n;\theta)$, the approximation $f_{n,m}$ converges at the same rate as $\cQ_{n} f$ as $n \rightarrow \infty$.  In addition, $m \geq \Theta(n;\theta)$ guarantees that $\| f_{n,m} \| \leq d_2 \theta^{-1} \| f \|$ and $\kappa(A) \leq d_2 \theta^{-1} \kappa(\tilde{A})$, thus making the linear system for $f_{n,m}$ solvable in a number of operations proportional to $\sqrt{d_2 \theta^{-1} \kappa(\tilde{A} )}$.

Note that $\Theta(n;\theta)$ is determined only by the sampling and reconstruction spaces $\rS_{m}$ and $\rT_{n}$.  Whilst $\Theta(n;\theta)$ can be numerically computed for any pair of spaces via the expression given in Lemma \ref{Cformlem}, analytical bounds must be determined on a case-by-case basis.  In the next section, where we consider the recovery of functions from their Fourier samples using (piecewise) polynomial bases, we are able to derive explicit forms for such bounds.  

\begin{remark}
As mentioned, the framework developed in this section was first introduced by the authors in \cite{BAACHShannon}.  Whilst a result similar to Theorem \ref{reconthm} was proved,  there are a number of important improvements offered by the theory presented in this paper:
\enum{
\item In \cite{BAACHShannon} it was assumed that the reconstruction vectors $\phi_1,\ldots,\phi_n$ were the first $n$ in an infinite sequence of vectors that formed a Riesz basis for $\rH$.  Conversely, Theorem \ref{reconthm} depends only on the subspace $\rT_{n}$, and thus the individual reconstruction vectors can be chosen arbitrarily.
\item The constants $K_{n,m}$ and $C_{n,m}$ are known exactly in terms of the sampling and reconstruction bases, and, whenever the sampling vectors are orthogonal, they possess a simple geometric interpretation in terms of subspace angles.  Moreover, these constants can be computed by determining either the minimal eigenvalue or, in certain cases, the norm of an $n \times n$ matrix.
\item Simple, explicit bounds for the condition number of the matrix $A$ are known in terms of the constant $C_{n,m}$ and the Gram matrix $\tilde A$.
\item The behaviour of $f_{n,m}$ as $m \rightarrow \infty$ (for $n$ fixed) can be fully explained in terms of oblique asymptotic optimality.
}
\end{remark}

\section{Polynomial reconstructions from Fourier samples}
\label{Fourrecovsec}
One of the most important examples of this procedure is the reconstruction of an analytic, but nonperiodic function $f$ to high accuracy from its Fourier coefficients.  Direct expansion in Fourier series converges only slowly in the $\rL^2(-1,1)$ norm, and suffers from the Gibbs phenomenon near the domain boundary.  Hence, given the first $m$ Fourier coefficients of $f$, we now seek to reconstruct $f$ to high accuracy in another basis.

Let $\rH = \rL^{2}(-1,1)$ (the space of square integrable functions on $(-1,1)$), $f : (-1,1) \rightarrow \bbR$ and
\bes{
\psi_{j}(x) = \frac{1}{\sqrt{2}} \E^{\I j \pi x},\quad j \in \bbZ,
}
be the standard Fourier basis functions.  For $m \geq 2$, we assume that the coefficients
\bes{
\hat{f}_{j} = \int^{1}_{-1} f(x) \overline{\psi_{j}(x)} \D x,\quad j = -\left \lfloor \frac{m}{2} \right \rfloor+1,\ldots,\left \lfloor \frac{m}{2} \right \rfloor-1,
}
are known (note that, whenever $m$ is even, this means that the first $m-1$ Fourier coefficients of $f$ are given.  We will allow this minor discrepancy since it simplifies ensuing analysis).  As a consequence of Theorem \ref{reconthm}, we are free to choose the reconstruction space.  The orthogonal projection of an analytic function onto the space $\bbP_{n-1}$ of polynomials of degree less than $n$ is known to converge exponentially fast at rate $\rho^{-n}$, where $\rho>1$ is determined from the largest Bernstein ellipse within which $f$ is analytic \cite{boyd}.  Hence, we let $\rT_{n} = \bbP_{n-1}$.  Note that an orthonormal basis for $\rT_{n}$ is given by the functions
\be{
\label{scaledLeg}
\phi_{j}(x) = \sqrt{j+\tfrac{1}{2}} P_{j}(x),\quad j \in \bbN,
}
where $P_{j}$ is the $j^{\rth}$ Legendre polynomial.  Moreover, if $\cQ_{n}$ is the orthogonal projection onto $\rT_{n}$, then it is well-known that
\be{
\label{SobErr}
\| f - \cQ_{n} f \| \leq c_{f} \sqrt{n} \rho^{-n},
}
where $c_f$ depends only on the maximal value of $f$ on the Bernstein ellipse indexed by $\rho$.  Naturally, we could also assume finite regularity of $f$ throughout, with suitable adjustments made to the various error estimates.  However, for simplicity we shall not do this.  

With this to hand, provided $m \geq \Theta(n;\theta)$ for some $\theta >0$, where $\Theta(n;\theta)$ 
is defined in \R{Theta}, the approximation $f_{n,m}$ obtained from the reconstruction procedure satisfies $\| f - f_{n,m} \| \leq c(\theta) \| f - \cQ_{n} f\|$ (see Theorem \ref{reconthm}).  In particular, $\| f - f_{n,m} \| \leq c(\theta) c_{f} \sqrt{n} \rho^{-n}$.  Hence, exponential convergence of $f_{n,m}$.  The key question remaining is how large $m$ must be in comparison to $n$ to ensure such rates of convergence.  Resolving this question involves estimating the quantity $C_{n,m}$, a task we next pursue.

\subsection{Estimates for $\Theta(n;\theta)$}
\label{Cnmest}
\begin{figure}
\begin{center}
$\begin{array}{ccc}
\includegraphics[width=6cm]{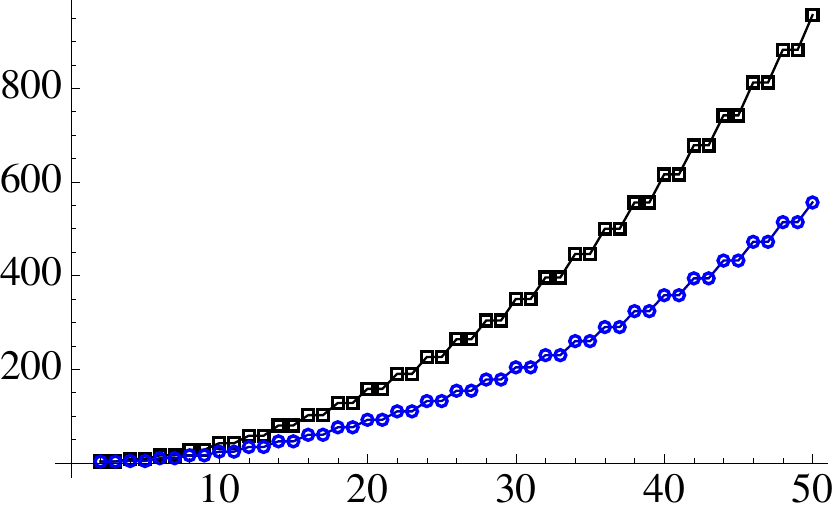} & \hspace{2pc} & \includegraphics[width=6cm]{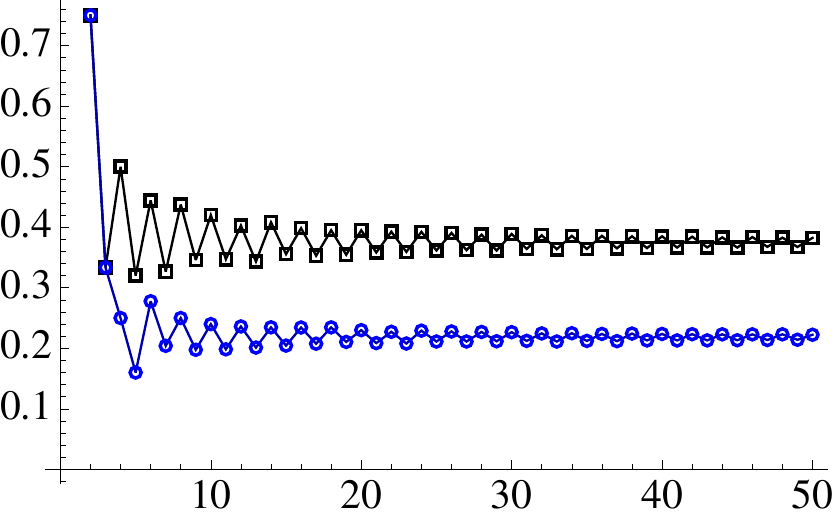}  
\end{array}$
\caption{\small The functions $\Theta(n;\theta)$ (left) and $n^{-2} \Theta(n;\theta)$ (right) for $\theta = \frac{1}{2}$ (squares) and $\theta = \frac{1}{4}$ (circles).} 
\label{LegPsiFig}
\end{center}
\end{figure}

For both numerical and analytical estimates of $\Theta(n;\theta)$ we need to select an appropriate basis of $\bbP_{n-1}$ (recall from Section \ref{genrecsect} that $\Theta(n;\theta)$ is independent of the basis used).  A natural choice is the orthonormal basis \R{scaledLeg} of scaled Legendre polynomials.  Fortunately, in this case, the inner products $\left < \phi_{k} , \psi_{j} \right >$ (i.e. the entries of the matrix $U$) are known in closed form:
\be{
\label{LegFC}
\left < \phi_{k} , \psi_{j} \right > = (-\I)^{k} \sqrt{\frac{k+\tfrac{1}{2}}{j}} J_{k+\frac{1}{2}}(j \pi),\quad j \in \bbZ,\ k \in \bbN,
}
where $J_{m}$ is the Bessel function of first kind.  This follows directly from
\be{
\label{BesseLegEq}
j_{m}(z) = \frac{1}{2} (-\I)^{m} \int^{1}_{-1} \E^{\I z x} P_{m}(x) \D x,\quad \forall z \in \bbC,
}
(see \cite[10.1.14]{AS}), where $j_{m}$ is the spherical Bessel function of the first kind, given by
\bes{
j_{m}(z) = \sqrt{\frac{\pi}{2 z}} J_{m+\frac{1}{2}}(z).
}
With this to hand, we may compute $C_{n,m}$ (and, in turn, $\Theta(n;\theta)$) via the expression given in Lemma \ref{Cformlem}.  In Figure \ref{LegPsiFig} we display the functions $\Theta(n;\frac{1}{2})$ and $\Theta(n;\frac{1}{4})$ against $n$.  Immediately, we witness quadratic growth of $\Theta(n;\theta)$ with $n$, a result we verify in this section.  In doing so, we shall also derive an upper bound for $\Theta(n;\theta)$ in terms of $n$ and $\theta$.  This gives an explicit, analytic condition for quasi-optimal recovery.  Whilst such a bound is completely robust (in that it holds for all $n$), we notice from Figure \ref{LegPsiFig} that $\Theta(n;\theta)$, when scaled by $n^{-2}$, quickly converges to an asymptotic limit.  In practice it is wasteful to use a larger value of $m$ than necessary (or, conversely, for fixed $m$ a overly pessimistic value of $n$).  Hence, in the second part of this section, we will also derive an asymptotic bound for $\Theta(n;\theta)$.

We commence as follows:
\lem{
\label{Cpolylem}
Suppose that $\rT_{n} = \bbP_{n-1}$, $\rS_{m}$ is the space spanned by the first $m$ Fourier basis functions and $m \geq \max \{2,\frac{2}{\pi} n\}$.  Then $C_{n,m}$ satisfies
\bes{
C_{n,m} \geq 1- \frac{4(\pi-2)n^2}{\pi^2 (2\lfloor \frac{m}{2} \rfloor-1)}. 
}
}
\prf{
From the definition of $C_{n,m}$ and the fact that $\{ \psi_{j} \}$ is an orthonormal basis we have
\bes{
1-C_{n,m} = 1 - \inf_{\substack{\phi \in \rT_{n} \\ \| \phi \| =1}} \left < \cP_{m} \phi , \phi \right > = \sup_{\substack{\phi \in \rT_{n} \\ \| \phi \| =1}} \left < \phi - \cP_{m} \phi , \phi \right >  = \sup_{\substack{\phi \in \rT_{n} \\ \| \phi \| =1}} \| \phi - \cP_{m} \phi \|^2,
}
where $\cP_{m}$ is the Fourier projection operator.  It now follows that $1-C_{n,m}\leq \sum^{n-1}_{k=0} \| \phi_k - \cP_{m} \phi_k \|^2$, where $\phi_k$ is given by \R{scaledLeg}.  By Parseval's theorem and the expression \R{LegFC}, we find that
\bes{
\| \phi_k - \cP_{m} \phi_k \|^2 = \sum_{|j| \geq \lfloor \frac{m}{2} \rfloor} \frac{k+\frac{1}{2}}{j} |J_{k+\frac{1}{2}}(j \pi) |^2.
}
Now, using a known result for Bessel functions \cite{hrycakIPRM}, it can be shown that 
\bes{
\frac{k+\frac{1}{2}}{j} |J_{k+\frac{1}{2}}(j \pi) |^2 \leq \frac{2k+1}{j \pi \sqrt{j^2 \pi^2-(k+\frac{1}{2})^2}},
}
provided $j \pi > k+\frac{1}{2}$.  Hence, for $m > \frac{2}{\pi} n$,
\be{
\label{phikerr}
\| \phi_k - \cP_{m} \phi_k \|^2 \leq \frac{2(2k+1)}{\pi^2} \sum_{j \geq \lfloor \frac{m}{2} \rfloor} \frac{1}{j \sqrt{j-\frac{(k+\frac{1}{2})^2}{\pi^2}}}.
}
Now, it was shown in \cite{hrycakIPRM} that $\sum_{j \geq m} \frac{1}{j \sqrt{j^2-c^2}} \leq \frac{1}{c} \arcsin \frac{c}{m-\frac{1}{2}}$, whenever $m \geq c +\frac{1}{2}$.  This gives
\bes{
\| \phi_k - \cP_{m} \phi_k \|^2 \leq \frac{4}{\pi} \arcsin \left [ \frac{2k+1}{(2 \lfloor \frac{m}{2} \rfloor-1) \pi} \right ],
}
and
\bes{
1-C_{n,m} \leq \frac{4}{\pi} \sum^{n-1}_{k=0} \arcsin \left [ \frac{2k+1}{(2 \lfloor \frac{m}{2} \rfloor-1) \pi} \right ].
}
We estimate this sum by the integral of $\arcsin t$.  We have
\bes{
1-C_{n,m} \leq 2(2 \lfloor \tfrac{m}{2} \rfloor-1) \int^{\frac{2n}{(2 \lfloor \frac{m}{2} \rfloor-1) \pi}}_{0} \arcsin t \D t .
}
Now, $\arcsin x \leq (\arcsin 1) x$ for $0 \leq x \leq 1$.  It follows that $F(x) = \int^{x}_{0} \arcsin t \D t \leq F(1) x^2$.  Computing this integral exactly, we arrive at $F(1) = \frac{\pi}{2}-1$.  Upon substituting $x = \frac{2n}{(2 \lfloor \frac{m}{2} \rfloor-1) \pi}$, this completes the proof.
}
This lemma confirms that it is sufficient for $m$ to scale quadratically with $n$ for quasi-optimal recovery.  Using this result, we find that
\thm{
\label{globbdthm}
Suppose that $\rT_{n}$ and $\rS_{m}$ are as in Lemma \ref{Cpolylem}.  Then, for $n \geq 2$, $\Theta(n;\theta)$ satisfies
\bes{
\Theta(n;\theta) \leq 2 \left \lceil \frac{1}{2}+\frac{2(\pi -2)}{\pi^2 (1-\theta)} n^2 \right \rceil,\quad \forall n \in \bbN.
}
}
\prf{
Suppose that $m \geq \{ 2 , \frac{2}{\pi} n \}$.  Then, by Lemma \ref{Cpolylem}, $C_{n,m} \geq \theta$ if
\bes{
1-\frac{4(\pi-2) n^2}{\pi^2 (2 \lfloor \frac{m}{2} \rfloor - 1)} \geq \theta.
}
Rearranging, we find that
\bes{
2 \left \lfloor \frac{m}{2} \right \rfloor \geq 1 + \frac{4(\pi-2) n^2}{\pi^2(1-\theta)}\quad \Rightarrow \quad m \geq 2 \left \lceil \frac{1}{2}+\frac{2(\pi -2)}{\pi^2 (1-\theta)} n^2 \right \rceil
}
and the theorem is proved, provided the right-hand side exceeds $\max \{ 2 , \frac{2}{\pi} \theta \}$.  Since $n \geq 2$, the right-hand side is certainly greater than $2$.  Moreover, 
\bes{
1+\frac{4(\pi-2) n^2}{\pi^2(1-\theta)} \geq \frac{8 (\pi-2) n}{\pi} > \frac{2n}{\pi},
}
as required.
}

\begin{figure}
\begin{center}
$\begin{array}{ccc}
\includegraphics[width=6cm]{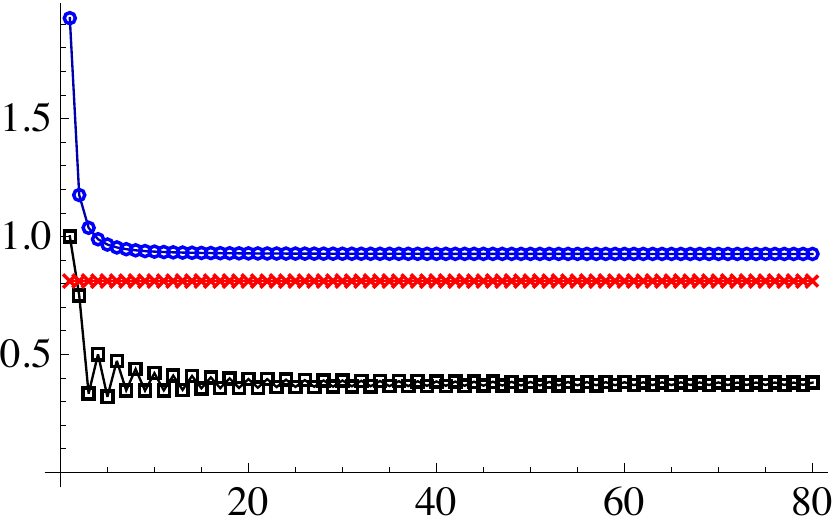} & \hspace{2pc} & \includegraphics[width=6cm]{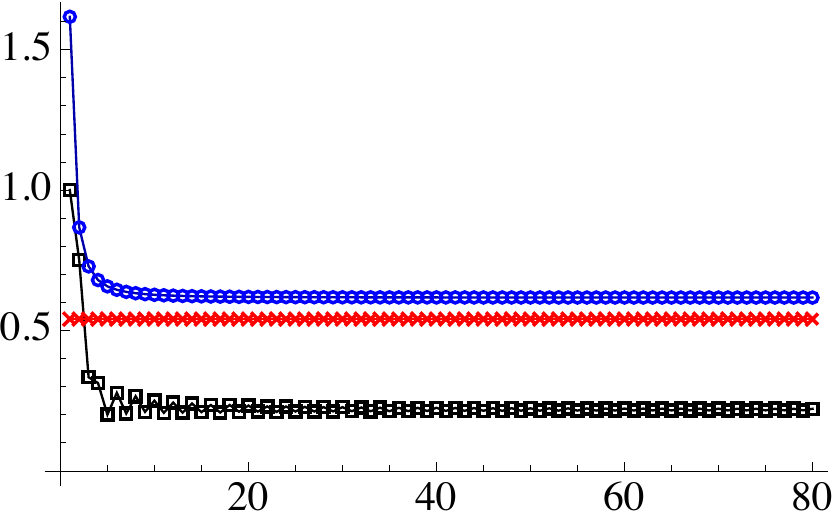}  
\end{array}$
\caption{\small The function $n^{-2} \Theta(n;\theta)$ (squares), the global bound (circles) and the asymptotic bound (crosses), for $n=2,\ldots,80$ and $\theta = \frac{1}{2}$ (left), $\theta = \frac{1}{4}$ (right).} 
\label{Phibdfig}
\end{center}
\end{figure}

Using a similar approach, we are also able to obtain an asymptotic bound for $\Theta(n;\theta)$ that is sharper than if were to use Theorem \ref{globbdthm} directly:

\thm{
\label{asybdthm}
Suppose that $\rT_{n}$ and $\rS_{m}$ are as in Lemma \ref{Cpolylem}.  Then the function $\Theta(n;\theta)$ satisfies
\bes{
n^{-2} \Theta(n;\theta) \leq \frac{4} {\pi^2 (1-\theta)} +\ord{n^{-2}},\quad n \rightarrow \infty.
}
}
\prf{
Suppose that $m = c n^2$ and recall \R{phikerr}.  Since $j<n$ and $k > \frac{1}{2} c n^2$, we deduce that
\bes{
\| \phi_k - \cP_{m} \phi_k \|^2 \leq \frac{2(2k+1)}{\pi^2} \sum_{j > \lfloor \frac{m}{2}  \rfloor} \frac{1}{j^2} + \ord{n^{-4}} = \frac{4(2k+1)}{c \pi^2 n^2} + \ord{n^{-4}}.
}
Hence
\bes{
1-C_{n,c n^2} \leq \frac{4}{c \pi^2 n^2} \sum^{n-1}_{k=0} (2k+1) + \ord{n^{-2}} = \frac{4}{c \pi^2 } + \ord{n^{-2}}.
}
Rearranging now gives the result.
}
In Figure \ref{Phibdfig} we compare the function $n^{-2} \Theta(n;\theta)$ for $\theta=\frac{1}{2},\frac{1}{4}$ and the global and asymptotic bounds of Theorems \ref{globbdthm} and \ref{asybdthm}.  Note that both bounds are reasonably sharp in comparison to the computed values.  In particular, as $n \rightarrow \infty$, $n^{-2} \Theta(n;\frac{1}{2})$ quickly approaches the limiting value $c \approx 0.38$, whereas the global and asymptotic upper bounds are $0.93$ and $0.81$ respectively.

At this moment, we reiterate an important point.  Whilst Legendre polynomials were used in the proof of Lemma \ref{Cpolylem}, the constant $C_{n,m}$ is independent of the particular reconstruction basis, and is only determined by the space $\rT_{n}$.  Hence, Theorems \ref{globbdthm} and \ref{asybdthm} provide \textit{a priori} estimates regardless of the particular implementation of the reconstruction procedure.  In the next section, we discuss the choice of polynomial basis and its effect on the numerical method.

\begin{remark}
\label{oversamplrmk}
In some applications, medical imaging, for example, oversampling is common.  Formally speaking, this is the situation where we wish to recover a function $f$ with support in $[-1,1]$ from its Fourier samples taken over an extended interval $K \supseteq [-1,1]$ (e.g. $K = [-\frac{1}{\epsilon},\frac{1}{\epsilon}]$ for some $0<\epsilon \leq 1$).  In this case, proceeding in a similar manner to before, we let $\rH = \rL^{2}(K)$,
\bes{
\psi_{j}(x) = \sqrt{\tfrac{c}{2}} \E^{\I c j \pi x},\quad x \in K,
}
where $c = \frac{1}{2} | K |$ and $\rT_{n} = \left \{ \phi : \phi |_{[-1,1]} \in \bbP_{n-1},\  \mbox{supp}(\phi) \subseteq [-1,1] \right \}$.  Using similar arguments to those of Lemma \ref{Cpolylem}, one can also derive estimates for $C_{n,m}$ and $\Theta(n;\theta)$ in this case.  In fact,
\be{
\label{Coversamp}
C_{n,m} \geq 1- \frac{4(\pi-2)n^2}{c \pi^2 (m-1)},
}
and
\bes{
\Theta(n;\theta) \leq 2 \left \lceil \frac{1}{2}+\frac{2(\pi -2)}{c \pi^2 (1-\theta)} n^2 \right \rceil,\ \forall n \in \bbN,\qquad n^{-2} \Theta(n;\theta) \leq \frac{4} {c \pi^2 (1-\theta)} +\ord{n^{-2}},\ n \rightarrow \infty.
}
In particular, we retain the scaling $m = \ord{n^2}$, regardless of the of size of the interval $K$.
\end{remark}

\subsection{Choice of polynomial basis}
\label{basischoicesect}
The results proved in this section are independent of the polynomial basis used for implementation.  In selecting such a basis, there are two questions which must be resolved.  First, how stable is the resultant method, and second, how can the entries of the matrix $U$ (as defined in Section \ref{comp_f}) be computed?  A straightforward choice is the orthogonal basis of Legendre polynomials \R{scaledLeg}.  In this case, $\tilde{A} = I$, where $\tilde{A}$ is the Gram matrix for $\{ \phi_0 , \ldots, \phi_{n-1} \}$, making the method well-conditioned (Lemma \ref{condnumlem}).  Moreover, the entries of $U$ are known explicitly via \R{LegFC}.  

Having said this, there is also interest in reconstructing in other polynomial bases.  In many circumstances it may be advantageous to have an approximation $f_{n,m}$ that is easily manipulable.  In this sense, an approximant composed of Legendre polynomials is not as convenient as one consisting of Chebyshev polynomials (of the first or second kind); the latter being easy to manipulate via the Fast Fourier Transform.  To this end, the purpose of this section is to detail the implementation of this method in terms of general Gegenbauer polynomials.

Gegenbauer polynomials arise as orthogonal polynomials with respect to the inner product
\bes{
\left < f , g \right >_{\lambda} = \int^{1}_{-1} f(x) \overline{g(x)} (1-x^2)^{\lambda-\frac{1}{2}} \D x,\quad  \lambda > -\frac{1}{2}.
}
For given $\lambda$, we denote the $j^{\rth}$ such polynomial by $C^{\lambda}_{j} \in \bbP_{j}$.  Important special cases are the Legendre polynomials ($\lambda =\frac{1}{2}$), and Chebyshev polynomials of the first ($\lambda = 0$) and second ($\lambda =1$) kind.  By convention \cite{bateman} (see also \cite{GottGibbs1}), each polynomial $C^{\lambda}_{j}$ is normalised so that
\be{
\label{Gegenendpt}
C^{\lambda}_{j}(1) = \frac{\Gamma(j+2 \lambda)}{j! \Gamma(2 \lambda)},
}
where $\Gamma$ is the Gamma function, in which case it is known that (see \cite[p.174]{bateman})
\be{
\label{Gegennorm}
\| C^{\lambda}_{j} \|^2_{\lambda} = \sqrt{\pi} \frac{\Gamma(j+2 \lambda) \Gamma(\lambda+\frac{1}{2})}{j! \Gamma(2 \lambda) \Gamma(\lambda) (j+\lambda)},
}
where $\| f \|_{\lambda} = \sqrt{\left < f , f \right >_{\lambda}}$.  With this to hand, we now define
\be{
\label{scaledGegen}
\phi_{j} = \frac{1}{\| C^{\lambda}_{j} \|_{\lambda}} C^{\lambda}_{j},\quad j=0,1,2,\ldots,
}
and seek to reconstruct $f$ in this basis.

Our first task is to compute the entries of the matrix $U$.  For this, we need to compute integrals of the form
\bes{
I_{k}(z) = \int^{1}_{-1} C^{\lambda}_{k}(x) \E^{\I z x} \D x, \qquad k= 0,1,2, \ldots,
}
where $z \in \bbR$.  Fortunately, such integrals obey a simple recurrence relation:
\lem{
\label{gegenrecurr}
For $z \neq 0$, the integrals $I_{k}(z)$ satisfy
\eas{
I_{0}(z) &= 2 C^{\lambda}_{0}(1) \frac{\sin z}{z},\quad I_{1}(z)=2 \I C^{\lambda}_{1}(1) \frac{\sin z - z \cos z}{z^2},
\\
I_{k+1}(z) &= \frac{2 \I(k+\lambda)}{z} I_{k}(z) + I_{k-1}(z) - \I\frac{ \E^{\I z}+(-1)^k \E^{-\I z} }{z} \left [ C^{\lambda}_{k+1}(1) - C^{\lambda}_{k-1}(1) \right ] , \quad k=1,2,\ldots.
}
When $z=0$, we have
\bes{
I_{0}(0) = 2 C^{\lambda}_{0}(1),\qquad I_{k}(0) = \frac{1 + (-1)^k }{2(k+\lambda)} \left [ C^{\lambda}_{k+1}(1) - C^{\lambda}_{k-1}(1) \right ] , \  k=1,2,\ldots.
}
}
\prf{
Recall the identity (see \cite[p.176]{bateman})
\bes{
C^{\lambda}_{j}(x) = \frac{1}{2(j+\lambda)} \frac{\D}{\D x} \left [ C^{\lambda}_{j+1} - C^{\lambda}_{j-1} \right ] ,\quad j=1,2,\ldots.
}
Substituting this into the expression for $I_k(z)$ and integrating by parts gives
\eas{
I_{k}(z) = \frac{1}{2(k+\lambda)} \left [ \left ( C^{\lambda}_{k+1}(x) - C^{\lambda}_{k-1}(x) \right ) \E^{\I z x}\right ]^{1}_{x=-1} - \frac{\I z}{2(k+\lambda)} \left [ I_{k+1}(z) - I_{k-1}(z) \right ].
}
Rearranging now gives the general recurrence for $k \geq 1$.  For $k=0,1$, we merely note that $C^{\lambda}_{0}(x) = C^{\lambda}_{0}(1)$, $C^{\lambda}_1(x) = C^{\lambda}_{1}(1) x$ and that
\bes{
\int^{1}_{-1} \E^{\I z x} \D x = 2 \frac{\sin z}{z},\quad \int^{1}_{-1} x \E^{\I z x} \D x = 2 \frac{\sin z - z \cos z}{z^2}.
}
The result for $z=0$ is derived in a similar manner.
}
Using this recurrence formula, the matrix $U$ can be formed in $\ord{m n}$ operations.  With this to hand, we now turn our attention to the condition number of $\tilde A$:
\thm{
\label{polyGramcond}
Let $\tilde A$ be Gram matrix for the vectors $\{ \phi_0,\ldots,\phi_{n-1} \}$, where $\phi_j$ is given by \R{scaledGegen}.  Then, $\kappa(\tilde A) = \ord{n^{|2\lambda-1|}}$ as $n \rightarrow \infty$.  In particular, whenever $\phi_0,\ldots,\phi_{n-1}$ arise from Chebyshev polynomials (of the first or second kinds), then $\kappa(\tilde A) = \ord{n}$. 
}
To prove this theorem, we first require the following two lemmas.  For convenience, we will write $\rL^{2}_{\lambda}(-1,1)$, $\lambda > -\frac{1}{2}$, for the space of square-integrable functions with respect to the Gegenbauer weight function $(1-x^2)^{\lambda-\frac{1}{2}}$.
\lem{
\label{weightednormlem}
Suppose that $-\frac{1}{2} < \lambda < \frac{1}{2}$.  Then, for all $g \in \rL^{\infty}(-1,1) $, we have $\| g \| \leq \| g \|_{\lambda}$ and, for some $c_{\lambda}>0$ independent of $g$,
\be{
\label{condneg}
\|  g\|_{\lambda} \leq c_{\lambda} \| g \|^{\lambda+\frac{1}{2}} \| g \|^{\frac{1}{2}-\lambda}_{\infty}.
}
Conversely, if $\lambda \geq \frac{1}{2}$ then, for all $g \in \rL^{\infty}(-1,1)$, $\| g \| \leq \| g \|_{\lambda} $ and
\be{
\label{condpos}
\| g \| \leq c_{\lambda}\| g \|^{\frac{1}{\lambda+\frac{1}{2}}}_{\lambda} \| g \|^{\frac{\lambda-\frac{1}{2}}{\lambda+\frac{1}{2}}}_{\infty}.
}
}
\prf{
Suppose first that $-\frac{1}{2} < \lambda < \frac{1}{2}$. Trivially, $\| g \| \leq \| g \|_{\lambda}$.  Now consider the other inequality.  For any $0<\epsilon<1$, we have
\eas{
\| g \|^2_{\lambda} &= \int^{1}_{-1} | g(x) |^2 (1-x^2)^{\lambda-\frac{1}{2}} \D x
\\
& = \int_{|x| \leq 1 - \epsilon} | g(x) |^2 (1-x^2)^{\lambda-\frac{1}{2}} \D x + \int_{ 1 - \epsilon<|x| \leq 1} | g(x) |^2 (1-x^2)^{\lambda-\frac{1}{2}} \D x
\\
& \leq (1-(1-\epsilon)^2)^{\lambda-\frac{1}{2}} \| g \|^2 +2 \| g \|^2_{\infty} \int^{1}_{1-\epsilon} (1-x^2)^{\lambda-\frac{1}{2}} \D x,
}
where $\| \cdot \|_{\infty}$ is the uniform norm on $[-1,1]$.  Note that $(1-(1-y)^2)^{\lambda-\frac{1}{2}} < y^{\lambda-\frac{1}{2}}$, $\forall y \in (0,1)$.  It now follows that
\bes{
\| g \|^2_{\lambda} \leq \epsilon^{\lambda-\frac{1}{2}} \| g \|^2 + \frac{2}{\lambda+\frac{1}{2}} \epsilon^{\lambda+\frac{1}{2}} \| g \|^2_{\infty} ,\quad 0 < \epsilon < 1.
}
Let $c > 2$ be arbitrary.  Then $\| g \|^2 < c \| g \|^2_{\infty}$, so we may let $\epsilon = \frac{\| g \|^2}{c\| g \|^2_{\infty}}$.  Substituting this into the previous expression immediately gives \R{condneg}.

Now suppose that $\lambda > \frac{1}{2}$.  Once more, trivial arguments give that $\| g \|_{\lambda} \leq \| g \|$.  For the other inequality, we proceed in a similar manner.  We have $\| g \|^2 \leq \epsilon^{\frac{1}{2}-\lambda} \| g \|^2_{\lambda}+ 2 \epsilon \| \phi \|^2_{\infty}$.  For $c>2$ we now set $\epsilon = \left ( \frac{\| g \|_{\lambda}}{c \| g \|_{\infty}} \right )^{\frac{2}{\lambda+\frac{1}{2}}}$, which gives \R{condpos}.
}

\lem{
\label{polynormineq}
Let $\lambda \geq \frac{1}{2}$.  Then, for all $\phi \in \bbP_{n-1}$, $\| \phi \|_{\infty}  \leq k_{n,\lambda} \| \phi \|_{\lambda}$, where $k_{n,\lambda}$ depends only on $n$ and $\lambda$ and satisfies $k_{n,\lambda} = \ord{n^{\lambda+\frac{1}{2}}}$ as $n \rightarrow \infty$.
}
\prf{
Let $\{ \phi_{0},\ldots,\phi_{n-1} \}$ be given by \R{scaledGegen}, and write an arbitrary $\phi \in \bbP_{n-1}$ as $\phi = \sum^{n-1}_{j=0} a_{j} \phi_{j}$, where $\sum^{n-1}_{j=0} | a_{j} |^2 = \| \phi \|^2_{\lambda}$.  By the Cauchy--Schwarz inequality,
\bes{
\| \phi \|^2_{\infty} \leq \| \phi \|^2_{\lambda}  \sum^{n-1}_{j=0}  \| \phi_j \|^2_{\infty} = k^2_{n,\lambda} \| \phi \|^2_{\lambda}. 
}
We now wish to estimate $\| \phi_{j} \|_{\infty}$.  Recall that $\| C^{\lambda}_{j} \|_{\infty}  = C^{\lambda}_{j}(1)$ \cite[p.206]{bateman}.  Hence, by \R{Gegenendpt} and \R{Gegennorm}, we have
\bes{
\| \phi_{j} \|^2_{\infty} = \frac{\Gamma(j+2\lambda) (j+\lambda)}{j! \sqrt{\pi} \Gamma(2 \lambda) \Gamma(\lambda+\frac{1}{2})}.
}
Consider the ratio $\frac{\Gamma(j+2 \lambda)}{j!}$.  By Stirling's formula,
\bes{
\frac{\Gamma(j+2 \lambda)}{j!} = \ord{j^{2\lambda-1}},\quad j \rightarrow \infty.
}
Hence $\| \phi_{j} \|^2_{\infty} = \ord{j^{2 \lambda}}$, which gives $k^2_{n,\lambda} = \ord{n^{2\lambda+1}}$, as required.
}

\prf{[Proof of Theorem \ref{polyGramcond}]
Since $\tilde A$ is Hermitian and positive definite, its condition number is the ratio of its maximum and minimum eigenvalues.  By a simple argument, we find that
\bes{
\lambda_{\max}(\tilde A) = \sup_{\substack{\phi \in \bbP_{n-1} \\ \phi \neq 0}} \frac{\| \phi \|^2}{\| \phi \|^2_{\lambda}},\quad  \lambda_{\min}(\tilde A) = \inf_{\substack{\phi \in \bbP_{n-1} \\ \phi \neq 0}} \frac{\| \phi \|^2}{\| \phi \|^2_{\lambda}}.
}
Consider the case $\lambda > \frac{1}{2}$.  By Lemma \ref{weightednormlem}, we have $\lambda_{\min}(\tilde{A}) \geq 1$ and
\bes{
\lambda_{\max}(\tilde A) \leq c^{2}_{\lambda} \sup_{\substack{\phi \in \bbP_{n-1} \\ \phi \neq 0}} \left ( \frac{\| \phi \|_{\infty}}{\| \phi \|_{\lambda}} \right )^{2 \frac{\lambda-\frac{1}{2}}{\lambda+\frac{1}{2}}}.
}
Using Lemma \ref{polynormineq}, we deduce that $\lambda_{\max}(\tilde A) = \ord{n^{2 \lambda - 1}}$, as required.  For the case $-\frac{1}{2} < \lambda < \frac{1}{2}$, we proceed in a similar manner.
}
This theorem confirms that the method can be implemented using Chebyshev polynomials whilst incurring only a mild growth in the condition number.  It follows that, if conjugate gradients are used to compute the approximation, the total computational cost of forming $f_{n,m}$ is $\cO(m n^{\frac{3}{2}})$, as opposed to $\ord{m n}$ in the Legendre polynomial case.  In the next section we present several examples of this implementation.

\begin{remark}
Whilst Theorem \ref{polyGramcond} provides an asymptotic estimate for $\kappa(\tilde A)$ (and hence $\kappa(A)$), it may also be useful to derive global bounds.  With effort, one could obtain versions of Lemmas \ref{weightednormlem} and \ref{polynormineq} involving explicit bounds.  For the sake of brevity, we shall not do this.  However, whenever Chebyshev polynomials are used (arguably the most important case), it is a relatively simple exercise to confirm that
\bes{
\kappa(\tilde A) \leq 2 \sqrt{2} n,\qquad \kappa(\tilde A) \leq 3^{\frac{2}{3}} \pi^{-\frac{1}{3}} \left [ n (n+\tfrac{1}{2}) (n+1) \right ]^{\frac{1}{3}},
}
in the first and second kind cases respectively.
\end{remark}

\subsection{Numerical examples}
\begin{figure}
\begin{center}
$\begin{array}{ccc}
\includegraphics[width=6cm]{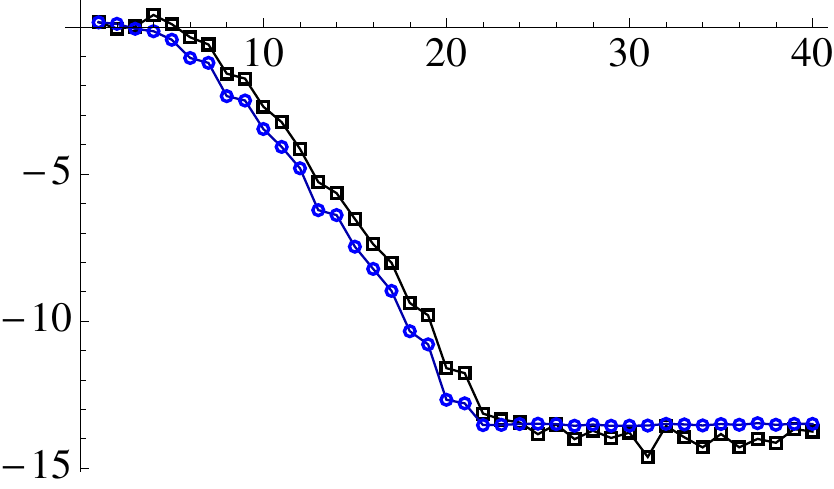} & \hspace{2pc} & \includegraphics[width=6cm]{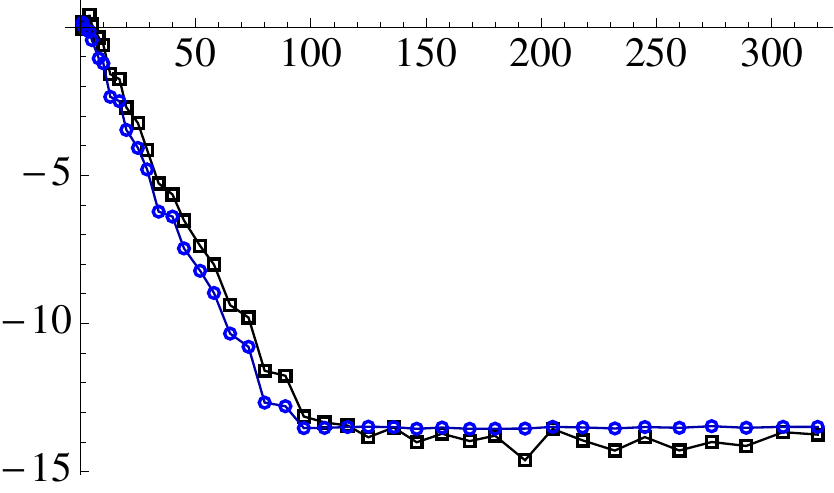}  
\end{array}$
\caption{\small Error in approximating $f(x) = \E^{-x} \cos 4x$ by $f_{n,m}(x)$ for $n=1,\ldots,40$.  Left: log error $\log_{10} \| f - f_{n,m} \|_{\infty}$ (squares) and $\log_{10} \| f - f_{n,m} \|$ (circles) against $n$.  Right: log error against $m = 0.2 n^2$.} 
\label{numexp1}
\end{center}
\end{figure}

We now present several numerical examples of this method.  All examples employ the value $m=0.2 n^2$, 
and the first series of examples consider the implementation using Legendre polynomials.  In Figure \ref{numexp1} we consider the function $f(x) = \E^{-x} \cos 4 x$.  Since $f$ is analytic in this case, we witness exponential convergence in terms of $n$ and root exponential convergence in terms of $m$.  Note the effectiveness of the method: using less than $100$ Fourier coefficients, we obtain an approximation with $13$ digits of accuracy.

\begin{figure}
\begin{center}
$\begin{array}{ccc}
\includegraphics[width=6cm]{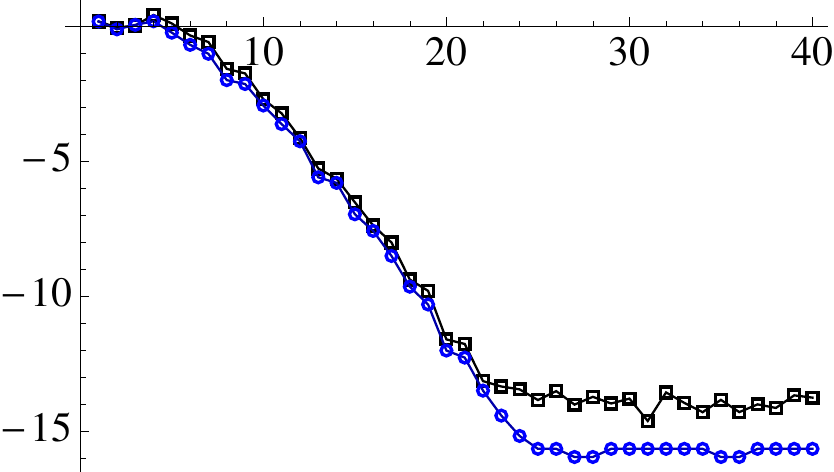} & \hspace{2pc} & \includegraphics[width=6cm]{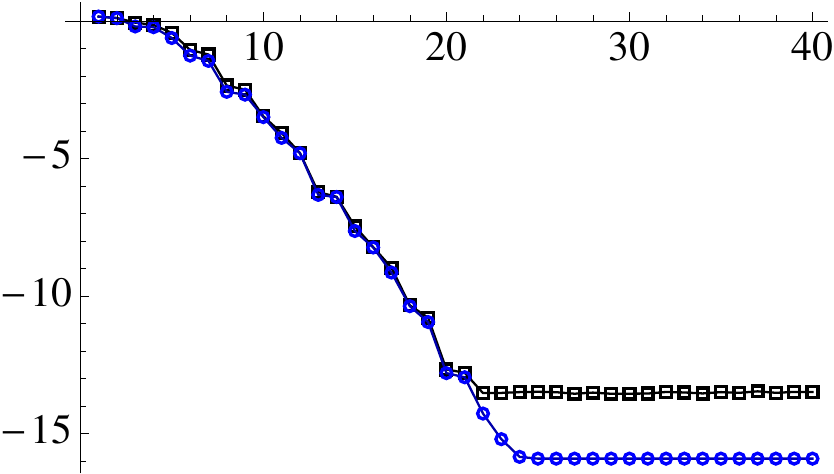}  
\end{array}$
\caption{\small Error in approximating $f(x) = \E^{-x} \cos 4x$ by $f_{n,m}(x)$ (squares) and $\cQ_{n}f(x)$ (circles) for $n=1,\ldots,40$.  Left: log uniform error.  Right log $\rL^2$ error.} 
\label{quasiopt1}
\end{center}
\end{figure}

As indicated by Theorem \ref{reconthm}, the approximation $f_{n,m}$ is quasi-optimal.  To highlight this feature of the method, Figure \ref{quasiopt1} displays both the error in approximating $f$ by $f_{n,m}$ and the best approximation $\cQ_{n} f$.  Note the very close correspondence of the two graphs.

\begin{figure}
\begin{center}
$\begin{array}{ccc}
\includegraphics[width=6cm]{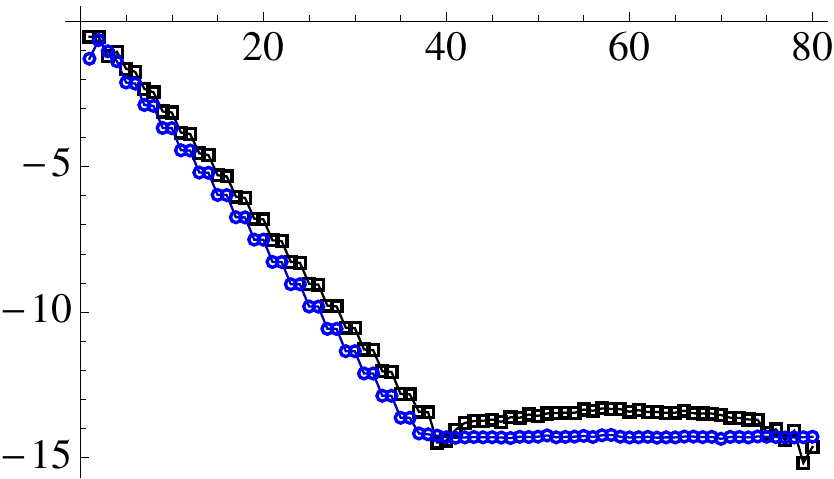} & \hspace{2pc} & \includegraphics[width=6cm]{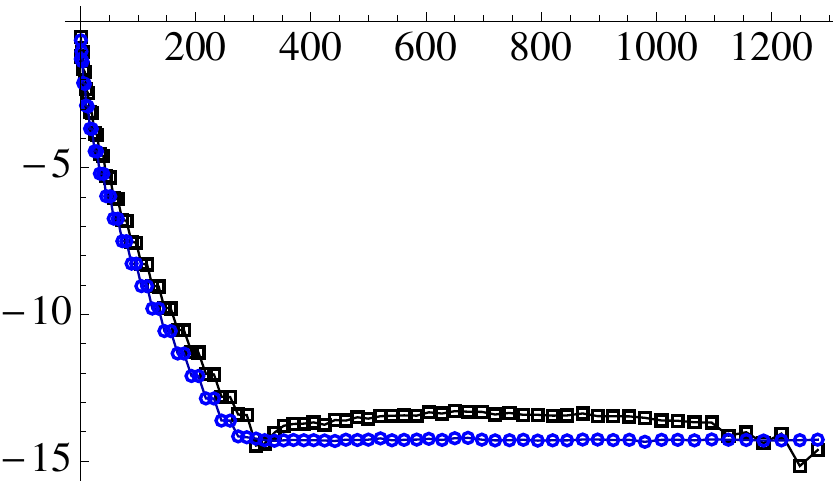}  
\end{array}$
\caption{\small Error in approximating $f(x) = \frac{1}{1+ x^2}$ by $f_{n,m}(x)$ for $n=1,\ldots,80$.  Left: log error $\log_{10} \| f - f_{n,m} \|_{\infty}$ (squares) and $\log_{10} \| f - f_{n,m} \|$ (circles) against $n$.  Right: log error against $m = 0.2 n^2$.} 
\label{numexp2}
\end{center}
\end{figure}

The example in Figures \ref{numexp1} and \ref{quasiopt1} is, in fact, entire.  Hence, the approximation $f_{n,m}$ converges super-geometrically in $n$ (as seen in Figure \ref{numexp1}).  For a meromorphic function, with complex singularity lying outside $[-1,1]$, the convergence rate is truly exponential at a rate $\rho$.  This is demonstrated in Figure \ref{numexp2}, the approximated function being $f(x) = \frac{1}{1+x^2}$.  Note that, despite the poles at $x=\pm \I$, the approximation $f_{n,m}$ still obtains $13$ digits of accuracy using only $250$ Fourier coefficients.

Next we consider reconstructions in other polynomials using the work of the previous section.  In Table \ref{ChebTab1} we give the error in approximating the function $f(x) = \E^{-x} \cos 4 x$ with Chebyshev polynomials of the first and second kinds.  Note that the resulting uniform error is virtually identical to the case of the Legendre polynomial implementation (unsurprisingly, since all three implementations compute exactly the same approximation $f_{n,m}$, up to numerical error).  Moreover, as evidenced by Table \ref{ChebTab2}, the payoff is only mild growth in the condition number $\kappa(A)$.

\begin{table}
\centering 
\begin{tabular}{c|ccccccccccc}
$n$ & 5 & 10 & 15 & 20 & 25 & 30 & 35 & 40   \\ \hline 
(a) & 1.45e0 & 1.85e-3 & 3.03e-7 &  2.53e-12 & 1.06e-14 & 8.42e-14 & 4.06e-14 & 5.31e-14\\
(b) & 1.45e0 & 1.85e-3 & 3.03e-7 & 2.53e-12 & 3.51e-14 & 1.16e-13 & 4.57e-14 & 7.70e-14\\ 
(c) & 1.45e0 & 1.85e-3 & 3.03e-7 & 2.49e-12 & 6.76e-14 & 7.33e-14 & 6.40e-14 & 5.15e-14 \\ 
\end{tabular} 
 \caption{\small Comparison of the error $\| f - f_{n,m} \|_{\infty}$ with $m=0.2 n^2$, where $f_{n,m}$ is formed from (a)  Legendre polynomials and Chebyshev polynomials of the (b) first and (c) second kinds.}
  \label{ChebTab1}
\end{table}

\begin{table}[t]
\centering 
\begin{tabular}{c|ccccccccccc}
$n$ & $5$ & $10$ & $15$ & $20$ & $25$ & $30$ & $35$ & $40$ \\ \hline 
(a) & $3.57$ & $5.55$ &$4.21$ & $5.20$ &$4.40$ & $5.06$ &$4.50$ & $6.77$\\
(b) &$13.74$ & $49.99$ & $52.63$ & $91.89$ & $92.89$ & $133.02$ & $133.49$ & $191.19$ \\ 
(c) & $3.90$ & $5.67$ & $7.25$ & $9.33$ & $11.91$ & $13.96$ & $16.56$  & $18.92$\\ 
\end{tabular}
 \caption{\small Comparison of the condition number $\kappa(A)$ with $m=0.2 n^2$, where $A$ is formed from (a)  Legendre polynomials and Chebyshev polynomials of the (b) first and (c) second kinds.}
  \label{ChebTab2}
\end{table}

\subsection{Connections to earlier work}
\label{earlierworksect}
Rather than choosing $m$ such that $C_{n,m} \geq \theta$, it may appear advantageous to find the minimum $m$ such that $C_{n,m} > 0$.  In other words, the smallest $m$ such that $f_{n,m}$ is guaranteed to exist.  Letting $\theta = 0$ in Theorems \ref{globbdthm} and \ref{asybdthm}, we immediately obtain a sufficient condition of the form $m \geq c n^2$, for some $c>0$.  However, this result is far too pessimistic: it is known that reconstruction is always possible, provided $m \geq n$ \cite{hrycakIPRM}.  For this reason, it may appear favourable to reconstruct using $m=n$, leading to the so-called \textit{inverse polynomial reconstruction method} \cite{shizgalGegen1,shizgalGegen2}.  Unfortunately, however, this approach is extremely unstable.  The linear system has geometrically large condition number, making the procedure extremely sensitive to both noise and round-off error.  Moreover, a continuous analogue of the Runge phenomenon occurs.  Roughly speaking, the approximation $f_{n,m}$ only converges to $f$ if geometric decay of $\| f - \cQ_{n} f \|$ is faster than the geometric growth of $\| A^{-1} \|$, meaning that only functions analytic in sufficiently large complex regions can be approximated by this procedure (as discussed in detail in \cite{BAACHShannon}, this behaviour can be understood in terms of the operator-theoretic properties of finite sections of certain non-Hermitian infinite matrices).  On the other hand, by allowing $m$ to range independently of $n$, we overcome all these difficulties, and obtain a stable method whose convergence is completely determined by the convergence of $\cQ_{n} f$ to $f$.

For the specific example of Legendre polynomial reconstructions from Fourier samples, this approach has also been recently considered in \cite{hrycakIPRM}.  Therein, the estimate $m = \ord{n^2}$ was derived, along with bounds for the error.  Naturally, this problem is just one specific example of our general framework.  However, within this context, our work improves and extends the results of \cite{hrycakIPRM}  in the following ways:
\enum{
\item Reconstruction is completely independent of the particular polynomial basis used.  In particular, the estimates for $\Theta(n;\theta)$ and $\| f - f_{n,m} \|$ are determined only by the spaces $\rT_{n}$ and $\rS_{m}$.  This allows for analysis of reconstructions in arbitrary polynomial bases, not just the
Legendre polynomials used in \cite{hrycakIPRM}.
\item The estimates for $\Theta(n;\theta)$ in Theorems \ref{globbdthm} and \ref{asybdthm} improve those given in \cite{hrycakIPRM}.  In particular, it was shown in \cite[Theorem 4.2]{hrycakIPRM} that
\be{
\label{hrycakest}
C_{n,\alpha n^2} \geq 1 - \frac{8}{\pi} \arcsin \frac{1}{\pi \alpha },\quad \forall n \in \bbN,\quad \alpha \geq 1,
}
(our constant $C_{n,m}$ corresponds to the quantity $\sigma^{2}_{n,m}$ in \cite{hrycakIPRM}).  Conversely, Lemma \ref{Cpolylem} leads to the improved bounds
\be{
\label{usglob}
C_{n, \alpha n^2} \geq 1 - \frac{4(\pi - 2)}{\pi^2 (\alpha - n^{-2})}\quad \forall n \geq \max \left \{ \frac{2}{\pi \alpha},\sqrt{\frac{2}{\alpha}} \right \},\quad \alpha > 0,
}
and
\be{
\label{usasy}
C_{n,\alpha n^2} \geq 1 - \frac{4}{\pi^2 \alpha} + \ord{n^{-2}},\quad n \rightarrow \infty,\quad \alpha > 0.
}
Not only are these bounds sharper, they also hold for a greater range of $\alpha$, thus permitting reconstruction with $m = \alpha n^2$ for any $\alpha > 0$, as opposed to just $\alpha \geq 1$.  This leads to savings in computational cost, and, in cases where $m$ is fixed, allows larger values of $n$ to be used, thereby increasing accuracy.  To illustrate this improvement, note that  \R{hrycakest} gives the estimate $\kappa(A) \leq 5.71$ when $m=n^2$.  Conversely, our estimate \R{usglob} yields the bound $2.61$ for $n \geq 2$, and \R{usasy} gives the asymptotic bound $1.68$.  To compare, direct computation of $\kappa(A)$ indicates that $\kappa(A) \leq 1.32$ for all $n$, and $\kappa(A) \rightarrow 1.2$ as $m \rightarrow \infty$.
\item Piecewise analytic functions and function of arbitrary numbers of variables can be recovered in a analogous fashion, with similar analysis (see Sections \ref{piecewisereconsect} and \ref{higherdimsect} respectively).
}

\subsection{Reconstruction of piecewise analytic functions}
\label{piecewisereconsect}

\begin{figure}
\begin{center}
$\begin{array}{ccc}
\includegraphics[width=6cm]{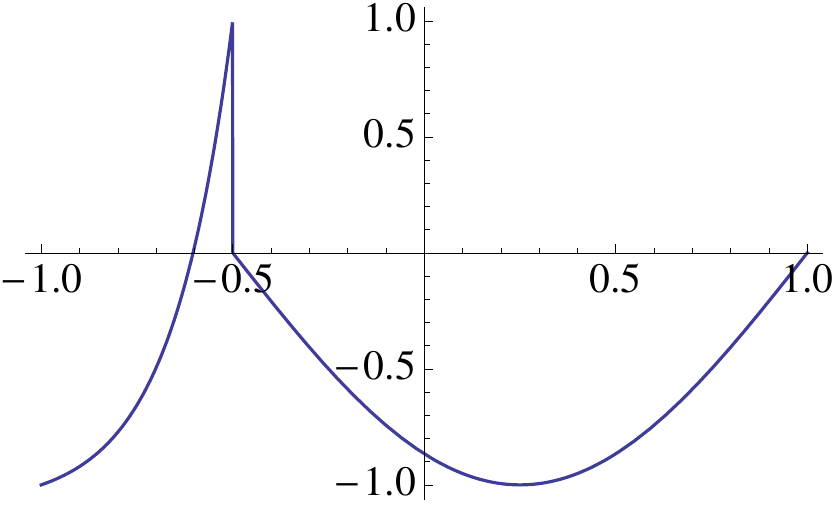} & \hspace{2pc} & \includegraphics[width=6cm]{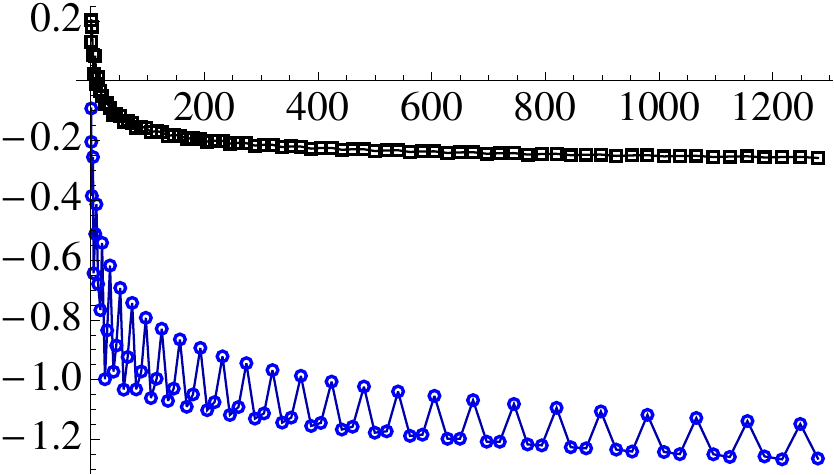}  
\end{array}$
\caption{\small Error in approximating the function \R{TannerTestFn} by $f_{n,m}(x)$ for $n=1,\ldots,80$.  Left: the function $f(x)$.  Right: log error $\log_{10} \| f - f_{n,m} \|_{\infty}$ (squares) and $\log_{10} \| f - f_{n,m} \|$ (circles) against $m = 0.2 n^2$.} 
\label{numexp5}
\end{center}
\end{figure}

Naturally, whenever the approximated function is not analytic, the convergence rate of the polynomial approximant $f_{n,m}$ to $f$ is not exponential.  For example, consider the function
\be{
\label{TannerTestFn}
f(x) = \left\{\begin{array}{lc} (2 \E^{2 \pi(x+1)} - 1 - \E^{\pi}) (\E^{\pi}-1)^{-1} & x \in [-1,-\frac{1}{2}) \\ - \sin (\frac{2 \pi x}{3} + \frac{\pi}{3}) &  x \in [-\frac{1}{2},1] \end{array}\right.
}
This function was put forth in \cite{tadmor2002adaptive} to test algorithms for overcoming the Gibbs phenomenon.  Aside from the discontinuity, its sharp peak makes it a challenging function to reconstruct accurately.  Since this function is discontinuous, we expect only low-order, algebraic convergence of $f_{n,m}$ in the $\rL^2$ norm, but no uniform convergence, an observation confirmed in Figure \ref{numexp5}.

However, by reconstructing this function in a polynomial basis, we are not exploiting the known information about $f$: namely, the jump discontinuity at $x=-\frac{1}{2}$.  The general procedure set out in Section \ref{genrecsect} allows us to use such information in designing a reconstruction basis.  Naturally, since $f$ is analytic in the subintervals $[-1,-\frac{1}{2}]$ and $[-\frac{1}{2},1]$, a better choice is to reconstruct $f$ in a piecewise polynomial basis.  The aim of this section is to describe this procedure.

Seeking generality, suppose that $f:[-1,1] \rightarrow \bbR$ is piecewise analytic with jump discontinuities at $-1<x_1<\ldots<x_l<1$.  Let $x_0=-1$ and $x_{l+1}=1$.  We assume that $f$ has been sampled via $\hat{f}_{j} = \left < f , \psi_{j} \right >$, $j=1,\ldots,m$, where $\left < \cdot , \cdot \right >$ is the Euclidean inner product on $\rL^{2}(-1,1)$.  In examples, these will be the Fourier samples of $f$, but the construction described below holds for arbitrary sampling bases consisting of functions defined on $[-1,1]$.

Throughout we shall assume that the discontinuity locations $x_1,\ldots,x_{l}$ are known exactly.  Although this may be a reasonable assumption in many applications, a fully-automated algorithm must also incorporate a scheme for singularity detection.  We shall not discuss possible approaches for doing this, and we refer the reader to \cite{Tadmor1} for further details.

Given the additional information about the location of the singularities of $f$, we now design a reconstruction basis to mirror this feature.  We shall construct such a basis via local co-ordinate mappings.  To this end, let $I_{r} = [x_r,x_{r+1}]$, $c_{r} = \frac{1}{2}(x_{r+1}-x_{r})$ and define $\Lambda_{r}(x) = \frac{x-x_r}{c_{r}}-1$, so that $\Lambda(I_r) = [-1,1]$.  Suppose now that $\rT'_{n}$ is a space of functions defined on $[-1,1]$ (e.g. the polynomial space $\bbP_{n-1}$).  By convention, we assume that each $\phi \in \rT'_{n}$ is extended by zero to the whole real line, i.e. $\phi(x) = 0$ for $x \in \bbR \backslash [-1,1]$.
 Let $\rT_{n,r}$ be the space of functions defined on $I_{r}$, given by $\rT_{n,r} =\left \{ \phi \circ \Lambda_{r} : \phi \in \rT'_{n} \right \}$.  We now define the new reconstruction space in the obvious manner:
\bes{
\rT_{n} = \left \{ \phi : \phi |_{I_r} \in \rT_{n_r,r},\  r=0,\ldots,l \right \},\quad n = \sum^{l}_{r=0} n_r,
}
and seek an approximation $f_{n,m} \in \rT_{n}$ to $f$ via the conditions $a_{m}(f_{n,m},\phi) = a_{m}(f,\phi)$, $\forall \phi \in \rT_{n}$, where $a_m$ is defined in (\ref{sesq}).   Suppose now that $\{ \phi_{1} , \ldots , \phi_{n} \}$ is a collection of linearly independent reconstruction functions with $\rT'_n = \mbox{span} \{ \phi_1 , \ldots, \phi_n \} $.  We construct a basis for $\rT_{n}$ by scaling.  To this end, we let $\phi_{r,j} = \frac{1}{\sqrt{c_r}}\phi_{j} \circ \Lambda_{r}$, and notice that $\rT_{n} = \mbox{span} \left \{ \phi_{r,j} : j=1,\ldots,n_r,\  r=0,\ldots,l \right \}$.  Note that, if $\{ \phi_j \}$ are orthonormal, then so are $\{ \phi_{r,j} \}$.  With this basis in hand, the approximation $f_{n,m}$ is now given by
\bes{
f_{n,m}= \sum^{l}_{r=0} \sum^{n_r}_{j=1} \alpha_{r,j} \phi_{r,j},
}
where the coefficients $\alpha_{r,j}$ are determined by the aforementioned equations.  As before, this is equivalent to the least squares problem $U \alpha \approx \hat{f}$ with block matrix $U = [ U_1 , \ldots, U_{l} ]$, where $U_{r}$ is the $m \times n_r$ matrix with $(j,k)^{\rth}$ entry
\bes{
\left < \phi_{r,k} , \psi_{j} \right > = \frac{1}{\sqrt{c_r}} \int^{x_{r+1}}_{x_r} \phi_{k}(\Lambda_{r}(x)) \psi_{j}(x) \D x.
\label{Uentry}
}
Here $\hat{f} = (\hat{f}_{1},\ldots,\hat{f}_{m})^{\top}$, $\alpha = [ \alpha_{0},\ldots,\alpha_{l} ] $ and $\alpha_{r} = (\alpha_{r,1},\ldots,\alpha_{r,n_r} )^{\top}$.  

Naturally, estimation of the constant $C_{n,m}$ is vital.  The following lemma aids in this task:
\lem{
\label{Clocalbdlem}
The constant $C_{n,m}$ satisfies
\bes{
C_{n,m} \geq 1 - \sum^{l}_{r=0} \left (1-C_{n_r,m,r} \right),
}
where $C_{n_r,m,r} = \inf_{\substack{\phi \in \rT_{n_r,r} \\ \| \phi \|=1 }} \left < \cP_{m} \phi , \phi \right >$.
}
\prf{
For $\phi \in \rT_{n}$, denote $\phi |_{I_{r}}$ by $\phi^{[r]}$.  Assume that $\phi^{[r]}$ is extended to $[-1,1]$ by zero, so that $\phi = \sum^{l}_{r=} \phi^{[r]}$.  Since $\phi^{[r]} \perp \phi^{[s]}$ for $r \neq s$, it follows that
\bes{
1-C_{n,m} = \sup \left \{ \frac{\sum^{l}_{r=0} \left <  \phi^{[r]} - \cP_{m} \phi^{[r]} , \phi^{[r]} \right >}{\sum^{l}_{r=0} \| \phi^{[r]} \|^2 } : \phi^{[r]} \in \rT_{n_r,r},\  r=0,\ldots,l,\ \sum^{l}_{r=0} \| \phi^{[r]} \|^2 \neq 0  \right \}. 
}
Note that, for $a_{r} \geq 0$ and $b_{r} > 0$, $r=0,\ldots,l$, the inequality
\bes{
\sum^{l}_{r=0} a_{r} \leq \sum^{l}_{r=0} \frac{a_r}{b_r} \sum^{l}_{r=0} b_{r},
}
holds.  Setting $a_{r} = \left <  \phi^{[r]} - \cP_{m} \phi^{[r]} , \phi^{[r]} \right >$ and $b_{r} = \| \phi^{[r]} \|^2$ and using this inequality gives
\eas{
1-C_{n,m} &\leq \sup \left \{ \sum^{l}_{r=0} \frac{\left <  \phi^{[r]} - \cP_{m} \phi^{[r]} , \phi^{[r]} \right >}{\| \phi^{[r]} \|^2 } : \phi^{[r]} \in \rT_{n_r,r},\  r=0,\ldots,l,\ \sum^{l}_{r=0} \| \phi^{[r]} \|^2 \neq 0  \right \}
\\
& \leq \sum^{l}_{r=0} \sup \left \{ \frac{\left < \phi - \cP_{m} \phi , \phi \right >}{\| \phi \|^2} : \phi \in \rT_{n_r,r},\  \| \phi \| \neq 0 \right \},
}
and this is precisely $ \sum^l_{r=0} (1-C_{n_r,m,r} )$.
}
Let us now focus on piecewise polynomial reconstructions from Fourier samples, in which case $\rT_{n}$ is the space of functions $\phi$ with $\phi |_{I_r}$ a polynomial of degree $n_{r}$.  Regarding the rate of convergence of the resulting approximation $f_{n,m}$, it is a simple exercise to confirm that
\bes{
\| f - f_{n,m} \| \leq c(\theta) c_f \sum^{l}_{r=0} \sqrt{n_{r}} \rho^{-n_{r}}_{r},
}
where $c(\theta)$ is defined in (\ref{smallc}), $c_f$ is a constant depending on $f$ only and $\rho_{r}$ is determined by the largest Bernstein ellipse (appropriately scaled) within which the function $f|_{I_r}$ is analytic.  
Hence, exponential convergence of $f_{n,m}$ to $f$.  The main question remaining is that of estimating the function $\Theta(n;\theta)$ for this reconstruction procedure.  For this, we have the following result, which extends Theorems \ref{globbdthm} and \ref{asybdthm} to this case:
\thm{
The function $\Theta(n;\theta)$ satisfies
\bes{
\Theta(n;\theta) \leq 2 \left \lceil \frac{1}{2}+\frac{2(\pi-2)}{\pi^2 (1-\theta)} \sum^{l}_{r=0} \frac{n^{2}_{r}}{c_r} \right \rceil,\quad \forall n = \sum^{l}_{r=0} n_{r},\ n_{0},\ldots,n_{l} \in \bbN,
}
and
\bes{
 \Theta(n;\theta) \leq \frac{4}{\pi^2 (1-\theta)} \sum^{l}_{0} \frac{n^2_r}{c_r} +\ord{1},\quad n_{0},\ldots,n_{l} \rightarrow \infty.
}
}
\prf{
In view of Lemma \ref{Clocalbdlem}, it suffices to consider $C_{n_r,m_r,r}$.  To this end, let $J = [\alpha , \beta ] \subseteq [-1,1]$, $\rT_{n,J}$ be the space of functions $\phi$ with $\mbox{supp}(\phi) \subseteq J$ and $\phi |_{J} \in \bbP_{n-1}$, and define
\bes{
E_{n,m} = \sup_{\substack{\phi \in \rT_{n,J} \\ \| \phi \|=1}} \left <  \phi-\cP_{m} \phi , \phi \right > .
}
Let $\Lambda(x) = \frac{x-\alpha}{c}-1$, where $c = \frac{1}{2} (\beta - \alpha)$, and write $\phi = \Phi \circ \Lambda$, where $\mbox{supp} (\Phi) \subseteq [-1,1]$.  Consider the quantity $\left < \phi , \psi_{j} \right >$.  By definition of $\psi_{j}$, we have
\bes{
\left < \phi , \psi_{j} \right > = \frac{1}{\sqrt{2}} \int^{1}_{-1} \phi(x) \E^{-\I j \pi x} \D x = \frac{c}{\sqrt{2} } \E^{-\I j \pi (\alpha + c)} \int^{\Lambda(1)}_{\Lambda(-1)} \Phi(y) \E^{-\I j \pi c y} \D y.
}
Let $K = [ \Lambda(-1) , \Lambda(1) ] = \Lambda([-1,1]) \supseteq [-1,1]$ and let $\cP_{m,K}$ be the Fourier projection operator based on the interval $K$. It now follows that
\bes{
E_{n,m} = \sup \left \{ \left < \Phi- \cP_{m,K} \Phi , \Phi \right > : \mbox{supp}(\Phi) \subseteq [-1,1], \  \Phi |_{[-1,1]} \in \bbP_{n-1},\  \| \Phi \| =1 \right \}.
}
This is now precisely the setup of Remark \ref{oversamplrmk}.  Using \R{Coversamp}, we therefore deduce that
\bes{
E_{n,m} \leq  \frac{4(\pi-2) n^2}{c \pi^2(2\lfloor \frac{m}{2} \rfloor-1)} .
}
Letting $J = I_{r}$, $c=c_r$ and using Lemma \ref{Clocalbdlem}, we now obtain
\be{
\label{pcwseCbd}
C_{n,m} \geq 1 - \frac{4(\pi-2)}{\pi^2(2\lfloor \frac{m}{2} \rfloor-1)} \sum^{l}_{r=0} \frac{n^2_r}{c_r},
}
from which the result follows immediately. 
}
To implement this scheme, it is necessary to compute the values \R{Uentry}.  By changing variables, it is easily seen that
\bes{
\left < \phi_{r,k} , \psi_{j} \right > = \sqrt{\frac{c_r}{2}} \E^{-\I j \pi d_r} \int^{1}_{-1} \phi_{k}(y) \E^{-\I j \pi c_r y} \D y,
}
where $d_r = \frac{1}{2}(x_{r+1}+x_{r})$.  Since \R{BesseLegEq} holds for all $z \in \bbC$, it follows that
\be{
\label{Upartcoeff}
\left < \phi_{r,k} , \psi_{j} \right > =  \E^{-\I j \pi d_r} (-\I)^{k} \sqrt{\frac{k+\frac{1}{2}}{j}} J_{k+\frac{1}{2}}(j \pi c_r),
}
whenever the functions $\phi_{r,k}$ arise from scaled Legendre polynomials.  Naturally, if the functions $\phi_{r,k}$ arise from arbitrary scaled Gegenbauer polynomials, computation of the values \R{Uentry} can be carried out recursively via the algorithm described in Section \ref{basischoicesect} (for appropriate choices of $z$).

\begin{figure}
\begin{center}
$\begin{array}{ccc}
\includegraphics[width=6cm]{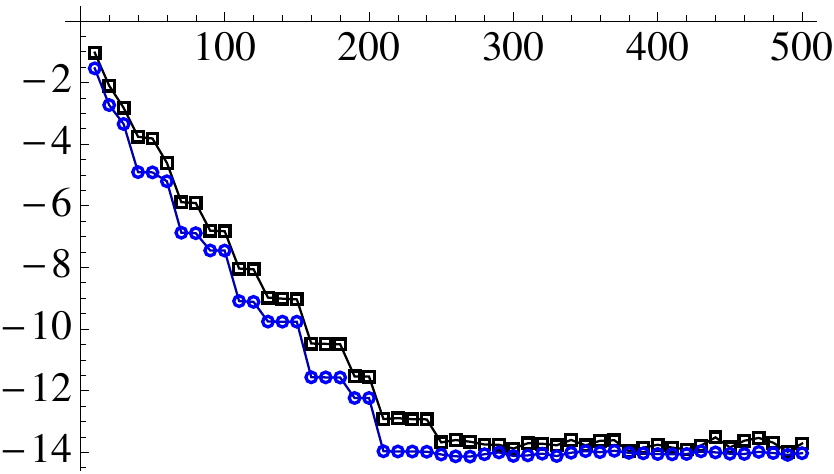} & \hspace{2pc} & \includegraphics[width=6cm]{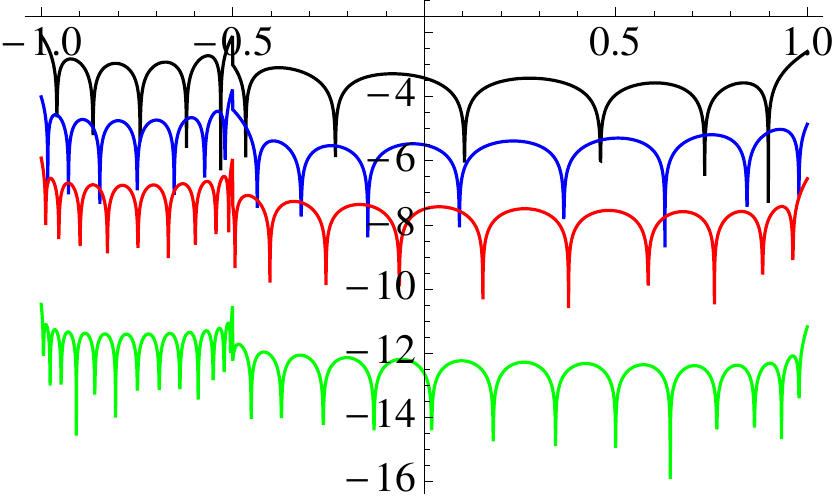}  
\end{array}$
\caption{\small Error in approximating the function \R{TannerTestFn} by $f_{n,m}(x)$.  Left: log error $\log_{10} \| f - f_{n,m} \|_{\infty}$ (squares) and $\log_{10} \| f - f_{n,m} \|$ (circles) against $m=1,\ldots,500$ with $n_0 = n_1$ chosen so that $m = \frac{1}{5} \left ( \frac{n^2_0}{c_0} + \frac{n^2_1}{c_1} \right )$.  Right: the error $\log_{10} |f(x)-f_{n,m}(x) |$ against $x \in [-1,1]$ for $m=20,40,80,160$.} 
\label{pcwsefig}
\end{center}
\end{figure}

In Figure \ref{pcwsefig} we apply this method to the function \R{TannerTestFn} using the orthonormal basis of scaled Legendre polynomials.  The improvement over Figure \ref{numexp5} is dramatic: using only $m \approx 250$ (with $n_0 = n_1 = 16$) we obtain $13$ digits of accuracy.  Note that, as expected, root exponential convergence occurs.  Moreover, as predicted by \R{pcwseCbd}, and illustrated in Table \ref{pcwsecond}, the condition number of the matrix $A$ remains small.

\begin{table}
\centering 
\begin{tabular}{c|ccccccccccc}
$m$ & $10$ & $20$ & $40$ & $80$ & $160$ & $320$ & $640$ & $1280$ \\ \hline 
$C_{n,m}$ & $0.05$ & $0.34$ &$0.33$ & $0.44$ &$0.44$ & $0.47$ &$0.49$ & $0.50$\\
$\kappa(A)$ & $19.88$ & $2.92$ &$3.06$ & $2.27$ &$2.27$ & $2.11$ &$2.03$ & $1.98$\\
\end{tabular}
 \caption{\small The constant $C_{n,m}$ and the condition number $\kappa(A)$ against $m$.}
  \label{pcwsecond}
\end{table}

\subsection{Comparison to existing methods}
\label{othermethsec}
Numerous methods exist for recovering functions to high accuracy from their Fourier data.  Applications are myriad, and range from medical imaging \cite{GelbMed2,GelbMed1} to postprocessing of numerical solutions of hyperbolic PDEs \cite{gottlieb2001spectral,jerri2007gibbs}.  Prominent examples which deliver high global accuracy (in contrast to standard filtering techniques, which only yield high accuracy away from the singularities of a function \cite{GottGibbsRev}) include Gegenbauer reconstruction \cite{GottGibbs3,GottGibbsRev,GottGibbs1}, techniques based on extrapolation \cite{Eckhoff3}, Pad\'e methods \cite{FourPade} and Fourier extension methods \cite{BoydFourCont,DHFEP} (for a more complete list, see \cite{boydlogsing} and references therein).

Whilst many of these methods deliver exponential convergence in terms of $m$ (the number of Fourier coefficients), they all suffer from ill-conditioning.  This comes as no surprise: the problem of reconstructing a function from its Fourier coefficients can be viewed as a continuous analogue of the recovery of a function from $m$ equidistant pointwise values.  As proved in \cite{TrefPlatteIllCond}, any method for this problem that converges exponentially fast in $m$ will suffer from exponentially poor conditioning.  We conjecture that a similar result holds in the continuous case.

Aside from increased susceptibility to round-off error and noise, ill-conditioning often makes so-called inverse methods (e.g. extrapolation and Fourier extension methods) costly to implement.  Conversely, the method proposed in the paper does not suffer from any ill-conditioning.  This negative consequence of \cite{TrefPlatteIllCond} is circumvented, precisely because we witness only root exponential convergence in $m$.  However, an advantage of this approach is that it delivers exponential convergence in $n$, the degree of the final approximant $f_{n,m}$.  In many applications it may be necessary manipulate $f_{n,m}$, its relatively low degree making such operations reasonably cheap.  Thus, this method has the advantage of compression, a feature not shared by the majority of the other methods mentioned previously.

A well-established and widely used alternative to this method is the Gegenbauer reconstruction technique, proposed by Gottlieb et al \cite{GottGibbs1}.  Much like this approach, it computes a polynomial approximation.  Yet it stands out as being direct, meaning that no linear system or least squares problem is required to be solved.  Whilst the original method has been shown to suffer from a generalised Runge phenomenon \cite{BoydGegen}, thereby severely affecting its applicability, an improved approach based on Freud polynomials was recently developed in \cite{GelbTannerGibbs}.

Comparatively speaking, this approach delivers exponential convergence in $\ord{m^2}$ operations.  On the other hand, our method obtains root exponential convergence at a cost of $\ordu{m^{\frac{3}{2}}}$ operations.  However, despite 
being theoretically more efficient, the various constants involved in the Gegenbauer/Freud procedure tend to be rather large.

Indeed, in Table \ref{CompTab} we compare the error in approximating the function \R{TannerTestFn} using both procedures (the data for the Freud method is taken from \cite[Table 1]{GelbTannerGibbs}.  Note that the parameter $N$ used therein is such that $2N = m$ is the total number of Fourier coefficients).  As is evident, the method proposed in this paper obtains an error of order $10^{-14}$ using less than $256$ Fourier coefficients, whereas the Freud procedure does not reach this value until more than $1024$ coefficients are used.

\begin{table}
\centering 
\begin{tabular}{c|ccccccc}
$m$ & 64 & 128 & 256 & 512 & 1024 & 2048 & 4096 \\ \hline 
(a) & 8.90e-01 & 1.37e-01 & 1.84e-04 & 1.01e-07 & 9.33e-13 & 5.27e-13 & 5.23e-14\\
(b) & 2.40e-04 & 8.36e-09 & 2.40e-14 & 1.38e-14 & 1.74e-14 & 2.26e-14 & 2.59e-14\ \\ 
\end{tabular}
 \caption{\small Comparison of the (a) Freud polynomial method and (b) generalised reconstruction method applied to \R{TannerTestFn}.  Here $m$ is the total number of Fourier samples.}
  \label{CompTab}
\end{table}

The most likely reason for this improvement is that the generalised reconstruction method is quasi-optimal, thereby delivering a near-optimal polynomial approximation, whereas the Gegenbauer and Freud procedures do not possess this feature.  Indeed, in the Gegenbauer case at least, we can quantitatively explain this phenomenon (the corresponding results for the Freud polynomial case are not so well understood).  It is known that, for analytic functions $f$, the approximation $F_{n,m}$ obtained from the Gegenbauer procedure converges geometrically at rate $\gamma^m$, where
\bes{
\gamma = \max \left \{ \rho^{-1} , \frac{(\beta + 2 \alpha)^{\beta + 2 \alpha}}{(2 \pi \E)^{\alpha} \alpha^{\alpha} \beta^{\beta}} \right \},
}
(see \cite[eqn. (4.12)]{GottGibbsRev}).  Here $\lambda = \alpha m$ is the parameter of the Gegenbauer polynomials used and $n = \beta m$ is the degree of the polynomial $F_{n,m}$.  Thus, in practice, the Gegenbauer method may converge much more slowly than the best polynomial approximation $\cQ_{m} f$, which converges at rate $\rho^{-m}$.  Conversely, our approach offers root exponential convergence at a rate determined solely by $\rho$: $\| f - f_{n,m} \| \sim \rho^{-\sqrt{m}}$.

Suppose now that the number of Fourier coefficients is fixed.  Then, ignoring various constants, the Gegenbauer method will offer a lower error than the generalised reconstruction procedure only when $\gamma^{m} \lesssim \rho^{-\sqrt{m}}$.  In other words,
\bes{
m \gtrsim \left ( \frac{\log \rho}{\log \gamma} \right )^2.
}
In typical examples (see \cite{GottGibbsRev,GottGibbs1}), the parameter values $\alpha = \beta = \frac{1}{4}$ have been used, giving the condition $m \gtrsim 19 (\log \rho)^2$.  Thus, for reasonably analytic functions $f$ (i.e. analytic in a sufficiently large Bernstein ellipse), the Gegenbauer procedure will only begin to outperform the generalised reconstruction method when the parameter $m$ is rather large.  Moreover, whenever $f$ is entire (as is the case with the example in Table \ref{CompTab}), the generalised reconstruction procedure will converge super-geometrically (in $n = \ord{\sqrt{m}}$), whereas the Gegenbauer method may still converge only exponentially at rate $\gamma$.  Thus, in this situation the Gegenbauer method may never outperform the generalised reconstruction method for any finite $m$.

Aside from improved numerical performance, let us mention several other benefits.  First, as discussed, the final approximation consists of only $\ordu{\sqrt{m}}$ terms, as opposed to $\ord{m}$.  Second, the basis for the polynomial reconstruction space $\rT_{n}$ can be chosen arbitrarily (in particular, independently of $m$) without affecting the convergence.  The only downside is a mild increase in condition number if nonorthogonal polynomials are employed.  In contrast, for the Freud/Gegenbauer procedure, only very specific types of polynomials can be used (which may not be simple to construct or manipulate \cite{GelbTannerGibbs}), and, whenever the number of samples $m$ is varied, all polynomials used for reconstruction must be changed.

To somewhat temper our claims, we mention that we have only considered one particular test function.  There may well be problems where the Freud/Gegenbauer procedure performs better, and the intention of future work is to present a more thorough comparison of the two methods.  That aside, one advantage of Gegenbauer method is that it is local: the approximation in each subdomain is computed separately and independently of any other subdomain.  Conversely, with our approach, the computations are inherently coupled.  Nevertheless, it may be possible to devise a local version of our approach, a question we intend to explore in future investigations.

\section{Reconstructions in tensor-product spaces}
\label{higherdimsect}
Thus far, we have focused on the reconstruction of univariate functions from Fourier samples.  A simple extension of this approach, via tensor products, is to functions defined in cubes.  The aim of this section is to detail this generalisation.

To formalise this idea, let us return to the general perspective of Section \ref{genrecsect}.  Suppose that the Hilbert space $\rH$ can be expressed as a tensor-product $\rH =  \rH_{1} \otimes \cdots \otimes \rH_{d}$ of Hilbert spaces $\rH_{i}$, $i=1,\ldots,d$, each having inner product $\left < \cdot , \cdot \right >_{i}$.  Note that, for $f = f_1 \otimes\cdots \otimes f_{d} \in \rH$ and $g = g_1 \otimes \cdots \otimes g_d \in \rH$, we have
\bes{
\left < f , g \right > = \prod^{d}_{i=1} \left < f_i , g_i \right >_i.
}
Now suppose that the sampling basis consists of tensor-product functions.  To this end, let
\bes{
\psi_{j} = \psi_{1,j_1} \otimes \cdots \otimes \psi_{d,j_d},\quad j = (j_1,\ldots,j_d) \in \bbN^d,
}
and, for $m = (m_1,\ldots,m_d) \in \bbN^d$, set
\bes{
\cS_{m} = \spn \left \{ \psi_{j} : j=(j_1,\ldots,j_d),\  1 \leq j_{i} \leq m_{i},\  i=1,\ldots, d \right \}.
}
We assume throughout that the collection $\{ \psi_{i,j} \}^{\infty}_{j=1}$ is a Riesz basis for $\rH_{i}$ for $i=1,\ldots,d$ (hence $\{ \psi_{j} \}$ is a Riesz basis for $\rH$).  With this to hand, we define the operator $\cP_{m} : \rH \rightarrow \cS_{m}$ by
\bes{
\cP_{m} f = \sum^{m_1}_{j_1=1} \cdots \sum^{m_d}_{j_d=1} \left < f , \psi_{j} \right > \psi_{j}.
}
Note that $\cP_{m} (f_1 \otimes \cdots \otimes f_d ) = \cP_{1,m_1} f_1 \otimes \cdots \otimes \cP_{d,m_d} f_{d}$, where $\cP_{i,m_i} : \rH_{i} \rightarrow \cS_{i,m_i}$ is defined in the obvious manner.  In a similar fashion, we introduce the reconstruction vectors $\phi_{j} = \phi_{1,j_1} \otimes \cdots \otimes \phi_{d,j_d}$, which form a basis for the reconstruction space
\bes{
\rT_{n} = \spn \left \{ \phi_{j} : j=(j_1,\ldots,j_d), \ 1 \leq j_i \leq n_i,\  i=1,\ldots,d \right \},\quad n = (n_1,\ldots,n_d) \in \bbN^d.
}
As before, we construct the approximation $f_{n,m} \in \rT_{n}$ via the equations $a_{m}(f_{n,m} , \phi ) = a_{m}(f,\phi)$, $\forall \phi \in \rT_{n}$, where $a_{m}(f,g) = \left < \cP_{m} f , g \right >$, $\forall f,g \in \rH$.

To cast this problem in a form suitable for computations, let $U^{[i]} \in \bbC^{m_i \times n_{i}}$ be the matrix with $(j,k)^{\rth}$ entry $\left < \psi_{i,j} , \phi_{i,k} \right >_i$.  Let $U \in \bbC^{\bar{m},\bar{n}}$ be the matrix of the $d$-variate reconstruction method, where $\bar{m} = m_1\ldots m_d$ and $\bar{n} = n_1 \ldots n_d$.  Then, it is easily shown that
\bes{
U = \bigotimes^{d}_{i=1} U^{[i]},\quad A = \bigotimes^{d}_{i=1} A^{[i]},
}
where $A = U^{\dag} U$, and $A^{[i]} = (U^{[i]})^{\dag} U^{[i]}$, and, in this case, $B_1\otimes B_2$ denotes the Kronecker product of matrices the $B_1$ and $B_2$.  By a trivial argument, we conclude that the number of operations required to compute $f_{n,m}$ is of order  $(n_1 m_1) \ldots (n_d m_d) \sqrt{ \kappa(A)}$.

Recall that the spectrum of the Kronecker product matrix $B_1\otimes B_2$ consists precisely of the pairs $\lambda \mu$, where $\lambda$ is an eigenvalue of $B_1$ and $\mu$ is an eigenvalue of $B_2$.  From this, we easily deduce that
\bes{
\kappa (A) = \prod^{d}_{i=1} \kappa  (A^{[i]}).
}
Hence, $\kappa(A)$ is completely determined by the matrices $A^{[i]}$, with the $i^{\rth}$ such matrix corresponding to the univariate reconstruction problem with sampling basis $\{ \psi_{i,j} \}^{m_i}_{j=1}$ and reconstruction basis $\{ \phi_{i,j} \}^{n_i}_{j=1}$.  Unsurprisingly, a similar observation also holds for the quantity $C_{n,m}$:
\lem{
\label{Ctensprodlem}
The constant $C_{n,m}$ satisfies $C_{n,m} = \prod^{d}_{i=1} C_{i,n_i,m_i}$.
}
\prf{
By Lemma \ref{Cformlem}, $C_{n,m} = \lambda_{\min}(\tilde{A}^{-1} A)$ and $C_{i,n_i,m_i} = \lambda_{\min}((\tilde{A}^{[i]})^{-1}A^{[i]})$, $i=1,\ldots,d$, where $\tilde{A}$ and $\tilde{A}^{[i]}$ are defined in the obvious manner.  Since $\tilde{A} = \tilde{A}^{[1]} \otimes \cdots \otimes \tilde{A}^{[d]}$, the matrix $\tilde{A}^{-1}$ is the Kronecker product of matrices $(\tilde{A}^{[i]})^{-1}$.  The result now follows immediately.
}

\subsection{Reconstruction of piecewise smooth functions}
Having presented the general case, we now turn our attention to the reconstruction of a piecewise smooth function $f:[-1,1]^d \rightarrow \bbR$.  We shall make the significant (see later discussion) assumption that $f$ is smooth in hyper-rectangular subregions of $[-1,1]^d$.  To this end, let $l_i \in \bbN$ for $i=1,\ldots,d$ and suppose that
\bes{
-1 = x_{0,i} < x_{1,i} < \ldots < x_{l_i,i} < x_{l_{i}+1,i} = 1,
}
and define $I_{r,i} = [x_{r,i},x_{r+1,i}]$ for $r=0,\ldots,l_{i}$.  For $r=(r_1,\ldots,r_d)$, we write $I_{r} = I_{r_1,1} \times \cdots \times I_{r_d,d}$, so that the collection
\bes{
 \left \{ I_{r} : r=(r_1,\ldots,r_d),\  r_i=0,\ldots,l_i,\  i=1,\ldots,d \right \},
}
consists of disjoint sets whose union is $[-1,1]^d$.  We assume that $f$ is smooth within each subdomain $I_r$.  In addition, for $r_i=0,\ldots,l_i$
 and $i=1,\ldots,d$, let $c_{r_i,i} = \frac{1}{2}(x_{r_i+1,i} -x_{r_i,i})$ and set $\Lambda_{r_i,i}(x) = \frac{x-x_{r_i,i}}{c_{r_i,i}}-1$, $x \in I_{r_i}$.  Note that $\Lambda_{r_i,i}(I_{r_i,i}) = [-1,1]$.

We now construct a reconstruction space.  To this end, for $n \in \bbN$ let $\rT'_n$, $\dim \rT'_{n}=n$, be a space of functions $\phi: \bbR \rightarrow \bbC$ with $\mbox{supp}(\phi) \subseteq [-1,1]$.   Define
\bes{
T_{n,r,i} = \left \{ \phi \circ \Lambda_{r,i} : \  \phi \in \rT'_n \right \},\quad n \in \bbN.
}
Now suppose that $n$ is the vector $(n_1,\ldots,n_d$), where
\bes{
n_i = \sum^{l_i}_{r=0} n_{r,i},\quad i=1,\ldots,d,
}
for some $n_{r,i} \in \bbN$.  We now define the reconstruction space $\rT_n$ by
\bes{
\rT_{n} = \bigotimes^{d}_{i=1}  \bigoplus^{l_i}_{r=0} T_{n_{r,i},r,i} .
}
We seek a basis for this space.  Let $\{ \phi_1,\ldots,\phi_n\}$, $n \in \bbN$, be a basis for $\rT'_n$, and set
\bes{
\phi_{r,j,i} = \frac{1}{\sqrt{c_{r,i}}} \phi_{j} \circ \Lambda_{r,i}.
}
A basis for $\rT_{n}$ is now given by
\bes{
\left \{ \phi_{r_1,j_1,1} \otimes \cdots \otimes \phi_{r_d,j_d,d},\  j=1,\ldots,n_{r_i,i},\  r_i=0,\ldots,l_i,\  i=1,\ldots,d \right \}.
}
This framework gives a general means in which to construct reconstruction bases in the tensor-product case for functions which are piecewise smooth with discontinuities parallel to the co-ordinate axes.  Suppose now that we consider the recovery of such a function from its Fourier samples.  Using the above framework, we are easily able to construct a basis consisting of piecewise polynomials of several variables.  The main question remaining is that of estimating the function $C_{n,m}$.  However, in view of the Lemma \ref{Ctensprodlem} and the results derived in Section \ref{Cnmest}, a simple argument gives

\begin{figure}
\begin{center}
$\begin{array}{ccc}
\includegraphics[width=6cm]{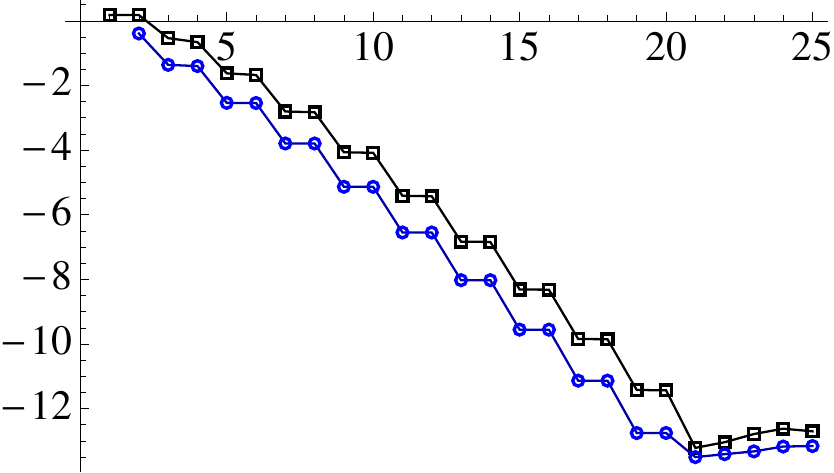} & \hspace{2pc} & \includegraphics[width=6cm]{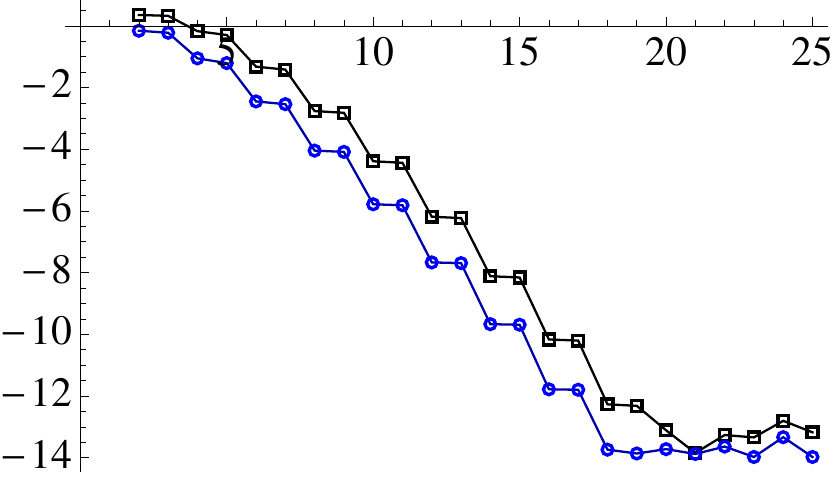}  
\end{array}$
\caption{\small The errors $\log_{10}\| f - f_{n,m} \|$ (squares) and $\log_{10} \| f - f_{n,m} \|_{\infty}$ (circles) for $n_{1} = n_{2}=1,\ldots,25$, where $f(x) = \E^{x^2 y}$ (left) and $f(x,y) = \sin 3 x y$ (right).} 
\label{Leg2DFig}
\end{center}
\end{figure}

\thm{
Suppose that $n =(n_1,\ldots,n_d)$, where $n_{i} = \sum^{l_i}_{r=0} n_{r,i}$, and let 
\bes{
\Theta(n;\theta) = \min \{ m = (m_1,\ldots,m_d) : C_{n,m} \geq \theta \},\quad \theta \in (0,1).
}
If $\theta_{1},\ldots,\theta_{d} \in (0,1)$ satisfy $\theta = \theta_1 \ldots \theta_d$, then we may write
\bes{
\Theta(n;\theta) = \left ( \Theta_{1}(n_1;\theta_1),\ldots,\Theta_{d}(n_d;\theta_d) \right ),
}
where, for $i=1,\ldots,d$, 
\bes{
\Theta_{i}(n_i ; \theta_i ) \leq 2 \left \lceil \frac{1}{2} + \frac{2(\pi-2)}{\pi^2(1-\theta_i)} \sum^{l_i}_{r=0} \frac{n^{2}_{r,i}}{c_{r,i}} \right \rceil,
}
and
\bes{
\Theta_{i}(n_i;\theta_i) \leq  \frac{4}{\pi^2(1-\theta_i)}\sum^{l_i}_{r=0} \frac{n^{2}_{r,i}}{c_{r,i}} +\ord{1},\quad n_{0,i},\ldots n_{l_i,i} \rightarrow \infty.
}
}
The main consequence of this theorem is the following: regardless of the dimension, the variables $m_{1},\ldots,m_d$ must scale quadratically with $n_{1},\ldots,n_{d}$ to ensure quasi-optimal recovery.  Consider now the most simple example of this approach: namely, where $f$ is smooth in $[-1,1]^d$, so that $\rT_{n}$ consists of multivariate polynomials. In Figure \ref{Leg2DFig} we plot the error in approximating the functions $f(x,y) = \E^{x^2 y}$ and $f(x,y) = \sin 3 x y$, using parameters $m_{1}=m_{2} = 0.5 n^2_{1}$ and $n_{2} = n_{1} $.  As in the univariate case, we observe the accuracy of this technique.  For example, using only $m_{1}=m_{2} \approx 200$ and $n_{1} = n_{2} \approx 20$ we obtain an error of order $10^{-14}$.

\begin{remark}
This approach (and many others based on tensor-product formulations) has the significant shortcoming that it requires the function to be singular in regions parallel to the co-ordinate axes.  Naturally, this is a rather restrictive condition.  For a function with singularities lying on a curve (in two dimensions, for example), one potential alternative is to apply the one-dimensional method along horizontal and vertical slices, and recover the two-dimensional function from the resulting one-dimensional reconstructions.

However, the generality of the reconstruction framework presented in this paper allows for another approach.  Given that we can reconstruct in any basis, an obvious alternative is to subdivide the domain into triangular elements and use a finite element basis for reconstruction.  This is a subject for future investigation.
\end{remark}

\section{Other sampling problems}
\label{othersampsect}
Overcoming the Gibbs phenomenon in Fourier series is an obvious application of the general framework developed in Section \ref{genrecsect}.  However, there is no reason to restrict to this case, and this framework can be readily applied to design effective methods for a variety of other problems.  In this section we describe several related problems, and the application of this framework therein.

\subsection{Modified Fourier sampling}
\label{modFoursamp}
Modified Fourier series were proposed in \cite{MF1} as a minor adjustment of Fourier series.  In the domain $[-1,1]$, rather than expanding a function $f$ in the classical Fourier basis
\bes{
\{ \cos j \pi x : j \in \bbN \} \cup \{ \sin j \pi x : j \in \bbN_{+} \},
}
we construct the modified Fourier expansion using the basis
\bes{
\{ \cos j \pi x : j \in \bbN \} \cup \{ \sin (j-\tfrac{1}{2}) \pi x : j \in \bbN_{+} \},
}
instead.  Though this basis arises from only a minor adjustment of the Fourier basis, the result is an improved approximation: the modified Fourier series of a smooth, nonperiodic function converges uniformly at a rate of $\ord{m^{-1}}$, whilst Fourier series suffers from the Gibbs phenomenon.  Although the convergence rate remains slow, the improvement over Fourier series, whilst retaining many of their principal benefits, has given rise to a number of applications of such expansions.  For a more detailed survey, we refer the reader to \cite{MFICOS}.

We shall consider modified Fourier expansions in a somewhat different context.  Given the similarity between the two bases, it is reasonable to assume that any sampling procedure (e.g. an MRI scanner) can be adapted to compute the modified Fourier coefficients of a given function (or image/signal), as opposed to the standard Fourier samples.  Indeed, if 
\bes{
\cF f(t) = \int^{1}_{-1} f(x) \E^{-\I \pi t x} \D x,
}
is the Fourier transform of $f$, then the modified Fourier coefficients are precisely
\eas{
\hat{f}^{C}_{j} &= \int^{1}_{-1} f(x) \cos j \pi x \D x =\frac{1}{2} \left [  \cF f(j) + \cF f(-j) \right ] ,
\\
\hat{f}^{S}_{j} &= \int^{1}_{-1} f(x) \sin (j-\tfrac{1}{2}) \pi x \D x =\frac{\I}{2} \left [  \cF f(j-\tfrac{1}{2}) + \cF f(\tfrac{1}{2}-j) \right ],
}
and hence can be computed from samples of the Fourier transform.  This raises the question: given that the general framework can handle sampling in either, is there an advantage gained from sampling in the modified Fourier basis, as opposed to the Fourier basis?  As we shall show, provided the function is analytic and nonperiodic, this is indeed the case.  Specifically, when we reconstruct in a polynomial basis, we require fewer samples to obtain quasi-optimal recovery to within a prescribed degree.

\begin{figure}
\begin{center}
$\begin{array}{ccc}
\includegraphics[width=6cm]{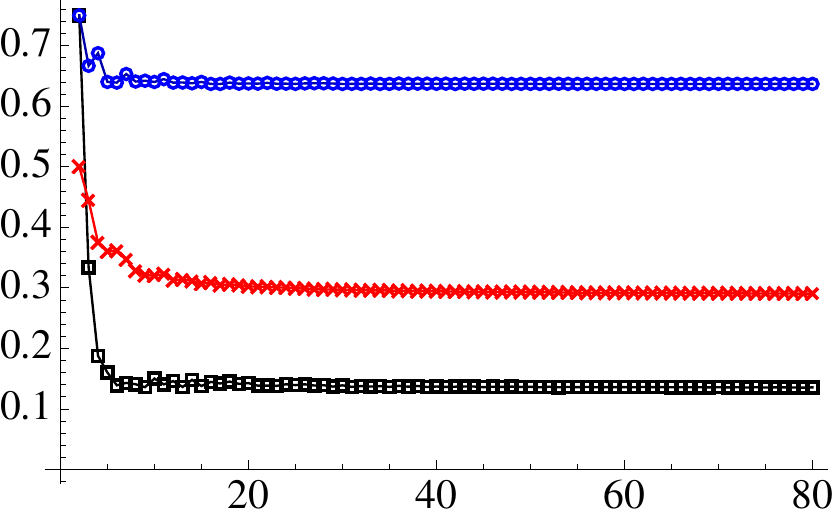} & \hspace{2pc} & \includegraphics[width=6cm]{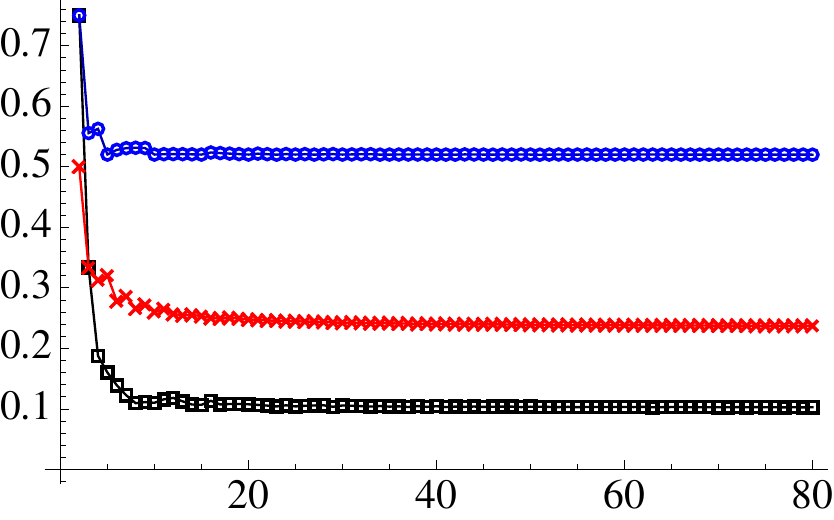}  
\end{array}$
\caption{\small The function $n^{-2} \Theta(n;\theta)$ (squares), the globabl bound (circles) and the asymptotic bound (crosses) for $n=1,\ldots,80$ and $\theta = \frac{1}{2}$ (left), $\theta = \frac{3}{4}$ (right).} 
\label{FMFcomp}
\end{center}
\end{figure}

Suppose that we carry out the reconstruction procedure as in Section \ref{Fourrecovsec} but using modified Fourier samples instead of Fourier samples.  For this, we set
\bes{
\cP_{m}f (x) = \frac{1}{2} \hat{f}^{C}_{0} + \sum^{\lfloor \frac{m}{2} \rfloor}_{j=1} \left [ \hat{f}^{C}_{j} \cos j \pi x + \hat{f}^{S}_{j} \sin j \pi x  \right ].
}
Naturally, we consider the function $\Theta(n;\theta)$ once more.  In Figure \ref{FMFcomp} we plot the function $\Theta(n;\theta)$ for the modified Fourier basis.  Upon comparison with Figure \ref{Phibdfig}, we conclude the following: using modified Fourier sampling, as opposed to standard Fourier sampling leads to a noticeable improvement.  Specifically, $n^{-2} \Theta(n ; \frac{1}{2})$ is approximately $0.15$ for large $n$ in the modified Fourier case, as opposed to $0.4$ in the Fourier case.

\begin{figure}
\begin{center}
$\begin{array}{ccc}
\includegraphics[width=6cm]{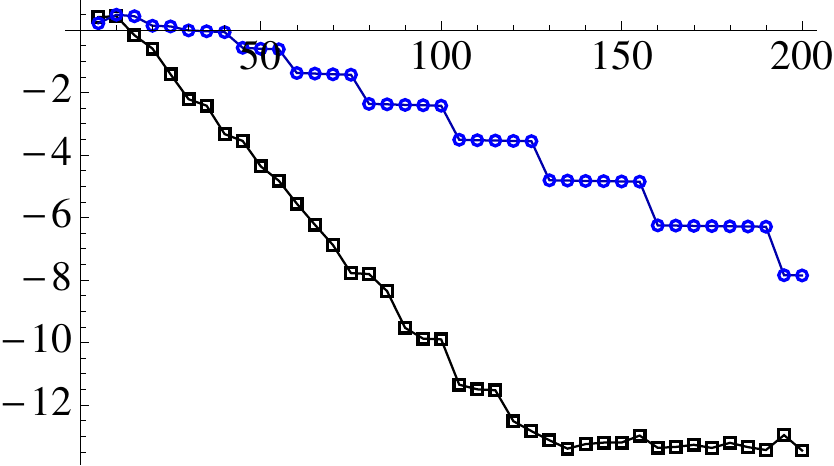} & \hspace{2pc} & \includegraphics[width=6cm]{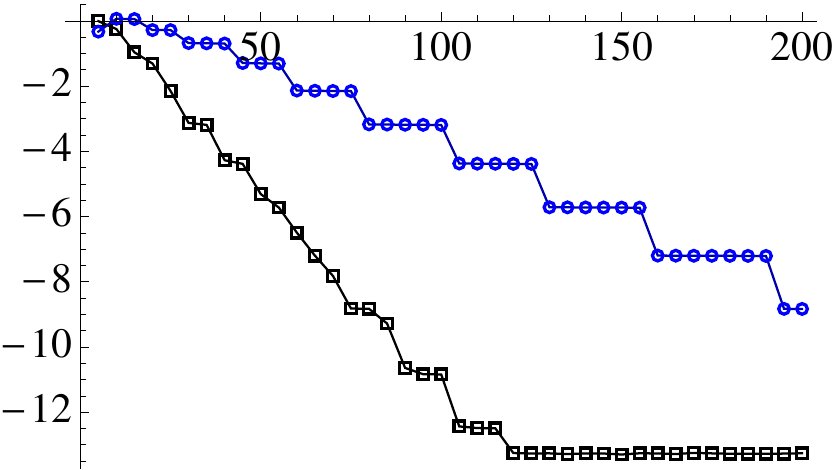}  
\end{array}$
\caption{\small The errors $\log_{10}\| f - f_{n,m} \|_{\infty}$ (left) and $\log_{10} \| f - f_{n,m} \|$ (right) against $m=5,10,15,\ldots,200$, where $f_{n,m}$ is computed from modified Fourier (squares) or Fourier (circles) samples.} 
\label{FMFcomp2}
\end{center}
\end{figure}

This result means that, if the number of samples $m$ is fixed, we are able to take a much larger value of $n$ in the modified Fourier case, whilst retaining quasi-optimal recovery (with constant $c(\theta)$).  As a result, we obtain better, guaranteed accuracy.  To illustrate this improvement, in Figure \ref{FMFcomp2} we compare the errors in approximating the function $f(x) = \E^{-x} \cos 8 x$ from either its Fourier or modified Fourier data.  In both cases the number of samples $m$ was fixed, and $n$ was chosen so that the parameter $C_{n,m} \geq \frac{1}{2}$.  As is evident, the method based on modified Fourier samples greatly outperforms the other.  For example, using only $m=120$ samples, we obtain an error of order $10^{-14}$ for the former, in comparison to only $10^{-4}$ for the latter.

As in the Fourier case, to implement the modified Fourier-based approach it is necessary to have estimates for the function $\Theta(n;\theta)$.  These are particularly simple to derive:
\lem{
\label{MFPhiBd}
Let $\rT_{n} = \bbP_{n-1}$ and $\rS_{m}$ be the space spanned by the first $m$ modified Fourier basis functions.  Then the function $\Theta(n;\theta)$ satisfies
\bes{
\Theta(n;\theta) \leq \left \lceil \frac{2}{\pi \sqrt{1-\theta}} k_{n} n^2 \right \rceil,\quad \mbox{where}\quad k_{n} = n^{-2} \sup_{\substack{\phi \in \bbP_{n-1} \\ \| \phi \|=1}} \| \phi' \|.
}
}
\prf{
In \cite{BA2} it was shown that $\| \phi - \cP_{m} \phi \| \leq \frac{2}{m \pi} \| \phi ' \|$ for all sufficiently smooth functions $\phi$.  The result now follows immediately from the definition of $C_{n,m}$.
}
As a result of this lemma, analytical bounds for $\Theta(n;\theta)$ are dependent solely on the constant $k_n$ of the Markov inequality $\| \phi' \| \leq k_n n^2 \| \phi \|$, $\forall \phi \in \bbP_{n-1}$.  Markov inequalities and their constants are well understood.  The question of determining $k_n$ was first studied by Schmidt \cite{schmidt1}, in which the estimates
 \be{
 \label{schmidtbds1}
k_{n} \leq \frac{1}{\sqrt{2}},\quad \forall n,\qquad \kappa_{n} \rightarrow \frac{1}{\pi},\quad n \rightarrow \infty,
}
were derived.  In \cite{schmidt2} the following improved asymptotic estimate was also obtained:
\be{
 \label{schmidtbds2}
k_{n} n^2 = \frac{(n+\frac{1}{2})^2}{\pi} \left [ 1-\frac{\pi^2-3}{12 (n+\frac{1}{2})^2} + \frac{R_{n}}{(n+\frac{1}{2})^4} \right ]^{-1},\quad n \geq 5,
}  
where  $-6<R_{n}<13$ (we refer the reader to \cite{BoettcherMarkov} for a more thorough discussion of both these results and more recent work on this topic).  Returning to $\Theta(n;\theta)$ we now substitute the result of Lemma \ref{MFPhiBd} into \R{schmidtbds1} and \R{schmidtbds2} to obtain the global and asymptotic bounds.  In Figure \ref{FMFcomp} we compare these bounds to their numerically computed values.  The relative sharpness of such estimates is once more observed.

\subsection{Polynomial sampling}
The primary concern of this paper has been reconstruction from Fourier (or Fourier-like) samples.  However, in several circumstances, most notably the spectral approximation of PDEs with discontinuous solutions \cite{GelbSGottChpt,gottlieb2001spectral}, the problem arises where a piecewise analytic function has been sampled in an orthogonal polynomial basis.  As previously noted, this approximation will converge slowly (and suffer from a Gibbs-type phenomenon near discontinuities), hence it is necessary to seek a new approximant with faster convergence.  Whilst a version of the Gegenbauer reconstruction method has been developed for this task \cite{GottGibbs4,GottGibbs3}, the advantages of the method proposed in this paper (see Section \ref{othermethsec}) make it a potential alternative to this existing approach.  Hence, the purpose of this section is to give a brief overview of this application.

It is beyond the scope of this paper to develop this example of the reconstruction procedure in its full generality.  Instead, we consider only the recovery of a piecewise analytic function $f:[-1,1] \rightarrow \bbR$ from its first $m$ Legendre polynomial coefficients $\hat{f}_{j} = \left < f , \psi_{j} \right >$, $j=0,\ldots,m-1$, where $\psi_{j} = (j+\frac{1}{2})^{\frac{1}{2}} P_{j}(x)$ is the $j^{\rth}$ normalised Legendre polynomial.  Proceeding as in Section \ref{piecewisereconsect}, we assume that $f$ has jump discontinuities as $-1<x_{1}< \ldots < x_{l} < 1$, and seek an approximation of the form
\bes{
f_{n,m}(x) = \sum^{l}_{r=0} \sum^{n_r-1}_{j=0} \alpha_{r,j} \phi_{r,j}(x),\quad n = \sum^{l}_{r=0} n_{r},
}
where $\phi_{r,j}(x) = \frac{1}{\sqrt{c_r}} \phi_{j}(\Lambda_{r}(x))$, $\Lambda_{r}(x) = \frac{x-x_r}{c_r}-1$, $c_{r} = \frac{1}{2}(x_{r+1}-x_r)$ and $\{ \phi_{0},\ldots,\phi_{n-1} \}$ is a system of polynomials on $[-1,1]$.  Since $f$ is piecewise analytic, we expect exponential convergence of $f_{n,m}$ to $f$, provided $m$ is sufficiently large in comparison to $n$.

\begin{figure}
\begin{center}
$\begin{array}{ccc}
\includegraphics[width=6cm]{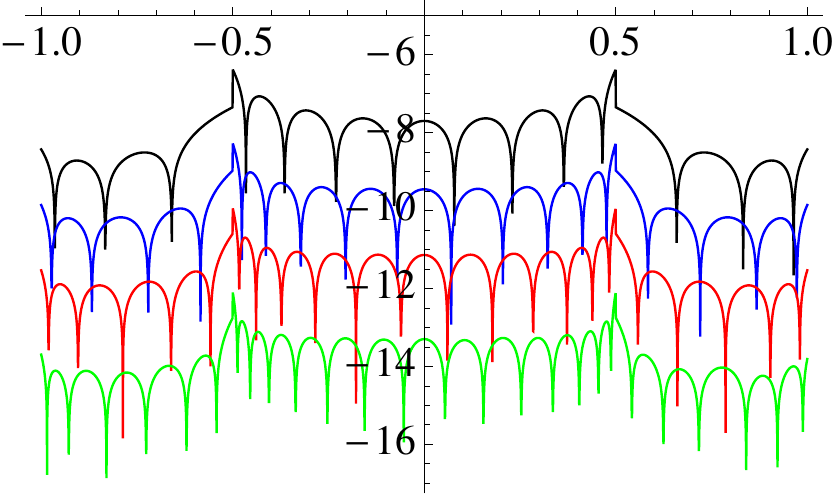} & \hspace{2pc} & \includegraphics[width=6cm]{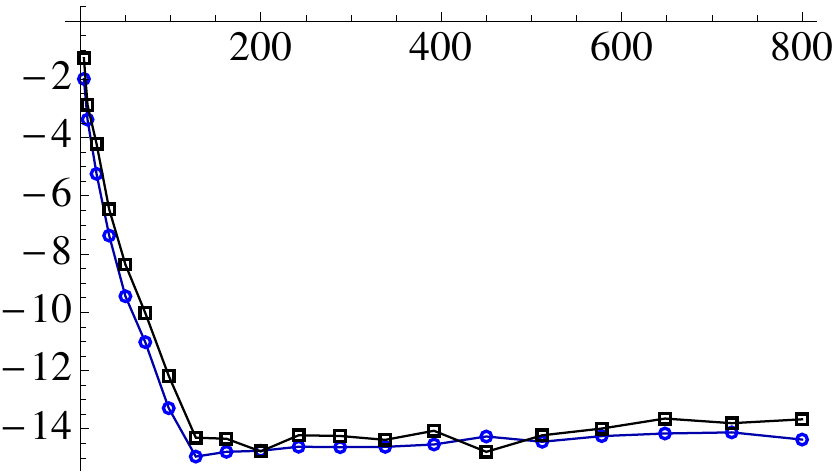}  
\end{array}$
\caption{\small Left: the error $\log_{10} |f(x)-f_{n,m}(x)|$ for $-1 \leq x \leq 1$ and $m=20,40,80,160$.  Right: log errors $\log_{10} \| f - f_{n,m} \|_{\infty}$ (squares) and $\log_{10} \| f - f_{n,m} \|$ (circles) against $m$.} 
\label{LegFig}
\end{center}
\end{figure}

Aside from determining how large $m$ must be in comparison to $n$ for recovery, the main question remaining is that of implementation, i.e. how to compute the entries of the matrix $U$.  This requires evaluation of the integrals
\bes{
\int^{x_{r+1}}_{x_r} \psi_{j}(x) \phi_{r,k}(x) \D x,\quad j=0,\ldots,m-1,\  k=0,\ldots,n-1.
}
Whenever the reconstruction functions $\phi_{r,k}$ arise from Gegenbauer polynomials, these calculations can be done iteratively.  For the sake of brevity, we will not describe this computation in full generality.  Instead, we consider only the situation where the functions $\phi_{r,k}$ arise from Legendre polynomials, in which case we are required to compute the integrals
\bes{
\int^{x_{r+1}}_{x_r} P_{j}(x) P_{k}(\Lambda_{r}(x)) \D x,\quad \quad j=0,\ldots,m-1,\  k=0,\ldots,n-1,\  r=0,\ldots,l.
}
We have
\lem{
\label{LegLegrecurr}
Let
\be{
\label{udef}
u_{j,k} = \int^{b}_{a} P_{j}(x) P_{k}(cx+d) \D x,\quad j,k=0,1,2,\ldots,
}
where $ca+d = -1$ and $cb+d = 1$.  Then
\bes{
u_{0,0} = b-a,\quad u_{j,0} = \frac{1}{j} \left [ P_{j+1}(x) - x P_{j}(x) \right ]^{b}_{x=a},\quad j=1,2,\ldots,
}
\bes{
u_{0,k} = \frac{2}{c} \delta_{0,k},\quad u_{1,k} = \frac{2}{3 c} \delta_{1,k} - \frac{2 d}{c} \delta_{0,k},\quad k=0,1,2,\ldots,
}
and, for $j \geq 2$ and $k \geq 1$,
\be{
\label{urecurr}
u_{j,k} = \frac{(2j-1)(k+1)}{c j (2k+1)} u_{j-1,k+1} + \frac{(2j-1) k}{c j (2k+1)} u_{j-1,k-1} - \frac{d(2j-1)}{c j} u_{j-1,k} - \frac{j-1}{j} u_{j-2,k}.
}
}
\prf{
Recall the recurrence relation 
\be{
\label{Legrec}
j P_{j}(x) = (2j-1) x P_{j-1}(x) - (j-1) P_{j-2}(x),\quad j=2,3,\ldots,
}
for Legendre polynomials \cite[chpt 22]{AS}.  Substituting this into \R{udef} gives
\bes{
u_{j,k} = \frac{2j-1}{j} \int^{b}_{a} x P_{j-1}(x) P_{k}(cx+d) \D x - \frac{j-1}{j} u_{j-2,k}.
}
Letting $x \mapsto cx+d$ in \R{Legrec} and rearranging, we find that
\bes{
x P_{k}(x) = \frac{k+1}{c(2k+1)} P_{k+1}(cx+d) - \frac{d}{c} P_{k}(cx+d) + \frac{k}{c(2k+1)} P_{k-1}(cx+d).
}
The recurrence \R{urecurr} now follows upon substituting this into the previous expression.

Next consider $u_{j,0}$.  Since $P_{0} \equiv 1$, we have $u_{0,0} = b-a$ and $u_{j,0} = \int^{b}_{a} P_{j}(x) \D x$ for $j \geq 1$.  Recall that the $j^{\rth}$ Legendre polynomial satisfies the Legendre differential equation \cite[chpt 22]{AS}
\bes{
\left [ (1-x^2) P'_{j}(x) \right ]' = - j(j+1) P_{j}(x).
}
Substituting for $P_{j}$ in $\int^{b}_{a} P_{j}(x) \D x$ and integrating gives
\bes{
u_{j,0} = \frac{1}{j(j+1)} \left [ (x^2-1) P'_{j}(x) \right ]^{b}_{x=a}.
}
The result now follows directly from the expression
\bes{
(1-x^2) P'_{j}(x) = (j+1)(x P_{j}(x) - P_{j+1}(x)),\quad j=0,1,2,\ldots.
}
To complete the proof, we consider $u_{0,k}$ and $u_{1,k}$.  By the assumptions on $a,b,c,d$, we find that
\bes{
u_{0,k} = \frac{1}{c} \int^{1}_{-1} P_{k}(x) \D x.
}
Orthogonality now gives $u_{0,k} = \frac{2}{c} \delta_{0,k}$, as required.  For $u_{1,k}$ we have
\bes{
u_{1,k} = \frac{1}{c} \int^{1}_{-1} (y-d) P_{k}(y) \D y = \frac{1}{c} \left ( \frac{2}{3} \delta_{1,k} - 2 d \delta_{0,k} \right ),
}
as required.
}
In Figure \ref{LegFig} we consider the approximation of the function
\be{
\label{Legtestfn}
f(x) = \left\{\begin{array}{ll} \sin \cos x & -\frac{1}{2} \leq x < \frac{1}{2}, \\  0 & \mbox{otherwise.} \end{array}\right.
}
by the aforementioned method, using parameter values $m = \frac{1}{8} n^2$, $n_{0} = n_{2} = \frac{1}{4} n$ and $n_1 = \frac{1}{2} n$.  As exhibited, we obtain $13$ digits of accuracy using only $m \approx 120$ Legendre coefficients of \R{Legtestfn}.  Note that, although we have not shown it, the scaling $m = \ord{n^2}$ appears to be sufficient for recovery.  Numerical results, demonstrating this hypothesis, are given in Table \ref{LegTab}.

\begin{table}
\centering 
\begin{tabular}{c|ccccccccccc}
$n$ & $8$ & $16$ & $24$ & $32$ & $40$ & $48$ & $56$ & $64$ & $72$ & $80$ \\ \hline 
$C_{n,m}$ & $0.98$ & $0.87$ & $0.85$ & $0.85$ & $0.85$ & $0.84$ & $0.84$ & $0.84$ & $0.84$ & $0.84$\\
$\kappa(U^{*} U)$ & $19.97$ & $4.17$ & $3.57$ & $3.57$ & $3.50$ & $3.43$ & $3.43$ & $3.41$ & $3.38$ & $3.38$\\ 
\end{tabular}
 \caption{\small The values $C_{n,m}$ and $\kappa(U^{*} U)$ against $n$, where $m = \frac{1}{8} n^2$.}
  \label{LegTab}
\end{table}

The function \R{Legtestfn} was introduced in \cite{GottGibbs4} to test the Gegenbauer reconstruction method when applied to this type of problem.  As shown in Figure \ref{LegFig}, we obtain a uniform error of roughly $10^{-8}$ using only $m = 40$ coefficients, and when $m=120$, the corresponding value is $10^{-14}$.  In comparison, the Gegenbauer method gives errors of roughly $10^{-3}$ and $10^{-7}$ for these values of $m$ (see \cite[Fig. 3]{GottGibbs4}), the latter being $10^{7}$ times larger. 

Whilst this method appears to be a promising alternative, it should be mentioned that the recursive scheme to compute the entries of $U$ requires $\ord{m^2}$ operations.  Since only $\ord{m n}$ operations are required to compute the approximation $f_{n,m}$, this is clearly less than ideal.  Having said that, the Gegenbauer reconstruction method requires $\ord{m^2}$ operations to compute each approximant, whereas with this scheme this higher cost is only incurred in a preprocessing step.

\section{Conclusions and future work}
\label{conclusions}
We have presented a reconstruction procedure to recover a function using any collection of linearly independent vectors, given a finite number of its samples with respect to an arbitrary Riesz basis.  This approach is both stable and quasi-optimal, provided the number of samples $m$ is sufficiently large in comparison to the number of reconstruction vectors $n$.  Moreover, this condition can be estimated numerically or, in certain circumstances, analytically.  A prominent example of this approach is the reconstruction of a piecewise analytic function from its Fourier samples.  Using a piecewise polynomial basis, this results in an approximation that converges root exponentially in terms of $m$, or exponentially in terms of $n$.

There are many potential avenues to pursue in extending this work, as we now detail:

\vspace{1pc}
\noindent \textit{1. Piecewise polynomial reconstructions from polynomial samples.}  In the final section in the paper we detailed the recovery of a piecewise analytic function in a piecewise polynomial basis, given its Legendre polynomial expansion coefficients.  Herein, an important open problem is verifying that the scaling $m = \ord{n^2}$ is sufficient for reconstruction.  Other challenges involve devising an iterative scheme for computing the entries of $U$ valid for reconstructions in arbitrary Gegenbauer polynomials, and which involves only $\ord{m n}$ operations, as opposed to $\ord{m^2}$.  Naturally, future work will also investigate the extension of this approach to reconstruction from arbitrary Gegenbauer polynomial expansion coefficients (as opposed to just Legendre polynomial expansion coefficients).

\vspace{1pc}
\noindent \textit{2. Gegenbauer polynomial reconstructions from Fourier samples.} As shown, the reconstruction procedure can be implemented with arbitrary Gegenbauer polynomials.  However, unless Legendre polynomials are used, the reconstruction is not stable.  This problem arises because Gegenbauer polynomials do not form a Riesz basis for the space $\rL^{2}(-1,1)$ unless $\lambda = \frac{1}{2}$.  However, Gegenbauer polynomial do form an orthogonal basis for the weighted space $\rL^{2}_{\omega}(-1,1)$, where $\omega(x) = (1-x^2)^{\lambda-\frac{1}{2}}$.  Hence, it is natural to ask whether the reconstruction procedure can be adjusted to incorporate this additional structure, thereby yielding a stable method.  It turns out that this can be done, with the first step being the derivation of an extended abstract framework along similar lines to Section \ref{genrecsect}.  We are currently compiling results in this case, and will report the details in a future paper.

\vspace{1pc}
\noindent \textit{3. Applications.}  Aside from the obvious applications in image and signal processing, the are a number of other potential uses of the procedure.  First, it may be applicable to the spectral discretisation of PDEs.  Spectral methods are extremely efficient for solving problems with smooth solutions.  However, for problems that develop discontinuities, e.g. hyperbolic conservation laws, a postprocessor is required to recover high accuracy \cite{gottlieb2001spectral}.  The Gegenbauer reconstruction technique has recently been applied to such problems (see \cite{GelbSGottChpt,gottlieb2001spectral} and references therein).  Given the potential advantages of the method developed in this paper (see Section \ref{othermethsec}), it is significant interest to apply this approach to these problems.  Aside from high accuracy, a pertinent issue in the use of spectral approximations for nonsmooth problems is the question of stability \cite{gottlieb2001spectral}.  Since the generalised reconstruction procedure is well-conditioned, there may also be a benefit in this regard.  Outside of PDEs, the Gegenbauer reconstruction technique has also been extended to other types of series, including radial basis functions \cite{jung2010recovery}, Fourier--Bessel series \cite{kamm-application} and spherical harmonics \cite{gelb1997resolution}.  Future work will also consider generalisation of the method of this paper along these lines. 

\vspace{1pc} \noindent
\textit{4. Spline and wavelet-based reconstructions.}
Reconstructions in spline and wavelet bases are important in numerous applications.  In \cite{BAACHShannon}, the authors gave a first insight into the application of such bases to the Fourier sampling problem.  However, the theory is far from complete.  Moreover, the use of more exotic objects, such a curvelets \cite{candes2004new} and ridgelets \cite{emmanuel1999ridgelets}, remains to be investigated.

\vspace{1pc} \noindent
\textit{5. Recovery from pointwise samples.}   The discrete analogue of the Fourier coefficient recovery problem involves the reconstruction of a function from $m$ equispaced samples in $[-1,1]$.  This problem has received more attention of late \cite{BoydRunge,TrefPlatteIllCond} than the continuous case considered in this paper.  In particular, the so-called mock--Chebyshev method \cite{boyd2009divergence} can be viewed as a discrete analogue of this approach.  Whilst the mock--Chebyshev method is well understood, there remain a number of other reconstruction from pointwise samples problems of interest.  In particular, with application to spectral collocation schemes based on Chebyshev or Legendre polynomials, the recovery of a piecewise analytic function from its values at Gauss or Gauss--Lobatto nodes.  Future work will consider this problem.

\small
\bibliography{AccRecovRefs}

\begin{thebibliography}{10}

\bibitem{AS}
M.~Abramowitz and I.~A. Stegun.
\newblock {\em Handbook of Mathematical Functions}.
\newblock Dover, 1974.

\bibitem{BA2}
B.~Adcock.
\newblock Multivariate modified {F}ourier series and application to boundary
  value problems.
\newblock {\em Numer. Math.}, 115(4):511--552, 2010.

\bibitem{BAACHShannon}
B.~Adcock and A.~C. Hansen.
\newblock A generalized sampling theorem for stable reconstructions in
  arbitrary bases.
\newblock {\em Technical report NA2010/07, DAMTP, University of Cambridge},
  2010.

\bibitem{MFICOS}
B.~Adcock and D.~Huybrechs.
\newblock Multivariate modified {F}ourier expansions.
\newblock In E.~R{\o}nquist et~al, editor, {\em Proceedings of the
  International Conference on Spectral and High Order Methods (to appear)},
  2010.

\bibitem{GelbMed2}
R.~Archibald, K.~Chen, A.~Gelb, and R.~Renault.
\newblock Improving tissue segmentation of human brain {M}{R}{I} through
  preprocessing by the {G}egenbauer reconstruction method.
\newblock {\em NeuroImage}, 20(1):489--502, 2003.

\bibitem{GelbMed1}
R.~Archibald and A.~Gelb.
\newblock A method to reduce the {G}ibbs ringing artifact in {M}{R}{I} scans
  while keeping tissue boundary integrity.
\newblock {\em IEEE Transactions on Medical Imaging}, 21(4):305--319, 2002.

\bibitem{bateman}
H.~Bateman.
\newblock {\em Higher Transcendental Functions}.
\newblock Vol. 2, McGraw--Hill, New York, 1953.

\bibitem{BoettcherMarkov}
A.~B{\"o}ttcher and P.~D{\"o}rfler.
\newblock Weighted {M}arkov-type inequalities, norms of {V}olterra operators,
  and zeros of {B}essel functions.
\newblock {\em Math. Nachr.}, 283(1):40--57, 2010.

\bibitem{boyd}
J.~P. Boyd.
\newblock {\em Chebyshev and Fourier Spectral Methods}.
\newblock Springer--Verlag, 1989.

\bibitem{BoydFourCont}
J.~P. Boyd.
\newblock A comparison of numerical algorithms for {F}ourier {E}xtension of the
  first, second, and third kinds.
\newblock {\em J. Comput. Phys.}, 178:118--160, 2002.

\bibitem{BoydGegen}
J.~P. Boyd.
\newblock Trouble with {G}egenbauer reconstruction for defeating {G}ibbs
  phenomenon: Runge phenomenon in the diagonal limit of {G}egenbauer polynomial
  approximations.
\newblock {\em J. Comput. Phys.}, 204(1):253--264, 2005.

\bibitem{boydlogsing}
J.~P. Boyd.
\newblock Acceleration of algebraically-converging {F}ourier series when the
  coefficients have series in powers of $1/n$.
\newblock {\em J. Comput. Phys.}, 228:1404--1411, 2009.

\bibitem{BoydRunge}
J.~P. Boyd and J.~R. Ong.
\newblock Exponentially-convergent strategies for defeating the {R}unge
  phenomenon for the approximation of non-periodic functions. {I}.
  {S}ingle-interval schemes.
\newblock {\em Commun. Comput. Phys.}, 5(2--4):484--497, 2009.

\bibitem{boyd2009divergence}
J.P. Boyd and F.~Xu.
\newblock Divergence ({R}unge phenomenon) for least-squares polynomial
  approximation on an equispaced grid and mock-{C}hebyshev subset
  interpolation.
\newblock {\em Appl. Math. Comput.}, 210(1):158--168, 2009.

\bibitem{emmanuel1999ridgelets}
E.~J. Cand{\`e}s and D.L. Donoho.
\newblock Ridgelets: a key to higher-dimensional intermittency?
\newblock {\em Phil. Trans. R. Soc. Lond. A}, 357(10):2495--2509, 1999.

\bibitem{candes2004new}
E.~J. Cand{\`e}s and D.L. Donoho.
\newblock New tight frames of curvelets and optimal representations of objects
  with piecewise ${C}^2$ singularities.
\newblock {\em Comm. Pure Appl. Math.}, 57(2):219--266, 2004.

\bibitem{christensen2003introduction}
O.~Christensen.
\newblock {\em An Introduction to Frames and {R}iesz Bases}.
\newblock Birkhauser, 2003.

\bibitem{FourPade}
T.~A. Driscoll and B.~Fornberg.
\newblock A {P}ad{\'e}-based algorithm for overcoming the {G}ibbs phenomenon.
\newblock {\em Numer. Algorithms}, 26:77--92, 2001.

\bibitem{Eckhoff3}
K.~S. Eckhoff.
\newblock On a high order numerical method for functions with singularities.
\newblock {\em Math. Comp.}, 67(223):1063--1087, 1998.

\bibitem{eldar2003sampling}
Y.C. Eldar.
\newblock Sampling without input constraints: Consistent reconstruction in
  arbitrary spaces.
\newblock {\em Sampling, Wavelets and Tomography}, 2003.

\bibitem{eldar2005general}
Y.C. Eldar and T.~Werther.
\newblock General framework for consistent sampling in {H}ilbert spaces.
\newblock {\em Int. J. Wavelets Multiresolut. Inf. Process.}, 3(3):347, 2005.

\bibitem{gelb1997resolution}
A.~Gelb.
\newblock The resolution of the {G}ibbs phenomenon for spherical harmonics.
\newblock {\em Math. Comp.}, 66(218):699--717, 1997.

\bibitem{GelbSGottChpt}
A.~Gelb and S.~Gottlieb.
\newblock The resolution of the {G}ibbs phenomenon for {F}ourier spectral
  methods.
\newblock In A.~Jerri, editor, {\em Advances in The Gibbs Phenomenon}. Sampling
  Publishing, Potsdam, New York, 2007.

\bibitem{GelbTannerGibbs}
A.~Gelb and J.~Tanner.
\newblock Robust reprojection methods for the resolution of the {G}ibbs
  phenomenon.
\newblock {\em Appl. Comput. Harmon. Anal.}, 20:3--25, 2006.

\bibitem{gilbarg2001elliptic}
D.~Gilbarg and N.S. Trudinger.
\newblock {\em Elliptic Partial Differential Equations of Second Order}.
\newblock Springer Verlag, 2001.

\bibitem{golub}
G.~H. Golub and C.~F. Van~Loan.
\newblock {\em Matrix Computations}.
\newblock John Hopkins University Press, Baltimore, 2nd edition, 1989.

\bibitem{gottlieb2001spectral}
D.~Gottlieb and J.~S. Hesthaven.
\newblock Spectral methods for hyperbolic problems.
\newblock {\em J. Comput. Appl. Math.}, 128(1-2):83--131, 2001.

\bibitem{GottGibbs4}
D.~Gottlieb and C-W. Shu.
\newblock On the {G}ibbs phenomenon {IV}: {R}ecovering exponential accuracy in
  a subinterval from a {G}egenbauer partial sum of a piecewise analytic
  function.
\newblock {\em Math. Comp.}, 64(211):1081--1095, 1995.

\bibitem{GottGibbs3}
D.~Gottlieb and C-W. Shu.
\newblock On the {G}ibbs phenomenon {III}: {R}ecovering exponential accuracy in
  a sub- interval from a spectral partial sum of a piecewise analytic function.
\newblock {\em SIAM J. Num. Anal.}, 33(1):280--290, 1996.

\bibitem{GottGibbsRev}
D.~Gottlieb and C-W. Shu.
\newblock On the {G}ibbs' phenomenon and its resolution.
\newblock {\em SIAM Rev.}, 39(4):644--668, 1997.

\bibitem{GottGibbs1}
D.~Gottlieb, C-W. Shu, A.~Solomonoff, and H.~Vandeven.
\newblock On the {G}ibbs phenomenon {I}: {R}ecovering exponential accuracy from
  the {F}ourier partial sum of a nonperiodic analytic function.
\newblock {\em J. Comput. Appl. Math.}, 43(1--2):91--98, 1992.

\bibitem{hrycakIPRM}
T.~Hrycak and K.~Gr\"{o}chenig.
\newblock Pseudospectral {F}ourier reconstruction with the modified inverse
  polynomial reconstruction method.
\newblock {\em J. Comput. Phys.}, 229(3):933--946, 2010.

\bibitem{DHFEP}
D.~Huybrechs.
\newblock On the {F}ourier extension of non-periodic functions.
\newblock {\em SIAM J. Numer. Anal.}, 47(6):4326--4355, 2010.

\bibitem{MF1}
A.~Iserles and S.~P. N{\o}rsett.
\newblock From high oscillation to rapid approximation {I}: Modified {F}ourier
  expansions.
\newblock {\em IMA J. Num. Anal.}, 28:862--887, 2008.

\bibitem{jerri1998gibbs}
A.J. Jerri, editor.
\newblock {\em The Gibbs phenomenon in Fourier Analysis, Splines, and Wavelet
  Approximations}.
\newblock Kluwer {A}cademic, {K}ordrecht, {T}he {N}etherlands, 1998.

\bibitem{jerri2007gibbs}
A.J. Jerri, editor.
\newblock {\em Advances in the Gibbs Phenomenon}.
\newblock Sampling Publishing, Potsdam, New York, 2007.

\bibitem{shizgalGegen1}
J.-H. Jung and B.~D. Shizgal.
\newblock Towards the resolution of the {G}ibbs phenomena.
\newblock {\em J. Comput. Appl. Math.}, 161(1):41--65, 2003.

\bibitem{shizgalGegen2}
J.-H. Jung and B.~D. Shizgal.
\newblock Generalization of the inverse polynomial reconstruction method in the
  resolution of the {G}ibbs phenomenon.
\newblock {\em J. Comput. Appl. Math.}, 172(1):131--151, 2004.

\bibitem{jung2010recovery}
J.H. Jung, S.~Gottlieb, S.O. Kim, C.L. Bresten, and D.~Higgs.
\newblock Recovery of high order accuracy in radial basis function
  approximations of discontinuous problems.
\newblock {\em J Sci Comput}, 45:359--381, 2010.

\bibitem{kamm-application}
J.R. Kamm, T.O. Williams, J.S. Brock, and S.~Li.
\newblock Application of {G}egenbauer polynomial expansions to mitigate {G}ibbs
  phenomenon in {F}ourier--{B}essel series solutions of a dynamic sphere
  problem.
\newblock {\em Int. J. Numer. Meth. Biomed. Engng.}, 26(1276--1292), 2010.

\bibitem{korner}
T.~W. K{\"o}rner.
\newblock {\em Fourier Analysis}.
\newblock Cambridge University Press, 1988.

\bibitem{TrefPlatteIllCond}
R.~Platte, L.~N. Trefethen, and A.~Kuijlaars.
\newblock Impossibility of fast stable approximation of analytic functions from
  equispaced samples.
\newblock {\em SIAM Rev. (to appear)}, 2010.

\bibitem{schmidt1}
E.~Schmidt.
\newblock Die asymptotische {B}estimmung des {M}aximums des {I}ntegrals
  {\"u}ber das {Q}uadrat der {A}bleitung eines normierten {P}olynoms, dessen
  {G}rad ins {U}nendliche w{\"a}chst.
\newblock {\em Sitzungsber. Preuss. Akad. Wiss.}, page 287, 1932.

\bibitem{schmidt2}
E.~Schmidt.
\newblock {\"U}ber die nebst ihren {A}bleitungen orthogonalen
  {P}olynomensysteme und das zugeh{\"o}rige {E}xtremum.
\newblock {\em Math. Ann.}, 119:165--204, 1944.

\bibitem{Tadmor1}
E~Tadmor.
\newblock Filters, mollifiers and the computation of the {G}ibbs' phenomenon.
\newblock {\em Acta Numerica}, 16:305--378, 2007.

\bibitem{tadmor2002adaptive}
E.~Tadmor and J.~Tanner.
\newblock Adaptive mollifiers for high resolution recovery of piecewise smooth
  data from its spectral information.
\newblock {\em Foundations of Computational Mathematics}, 2(2):155--189, 2002.

\bibitem{unser2000sampling}
M.~Unser.
\newblock Sampling--50 years after {S}hannon.
\newblock {\em Proc. IEEE}, 88(4):569--587, 2000.

\end{thebibliography}

\end{document}